\documentclass[11pt]{amsart}  
\usepackage[T1]{fontenc}
\usepackage{amsxtra}

\usepackage{color}
\usepackage{amscd,amsfonts,amsmath,amssymb,
amsthm,amstext,latexsym}
\usepackage{enumerate,pifont,graphicx}
\usepackage[english]{babel}
\usepackage[all,cmtip]{xy} 
\newcommand{\A}{{\mathbb A}}
\newcommand{\Z}{{\mathbb Z}}
\renewcommand{\S}{{\mathbb S}}
\newcommand{\R}{{\mathbb R}}
\newcommand{\C}{{\mathbb C}}
\newcommand{\D}{{\mathbb D}}
\newcommand{\N}{{\mathbb N}}
\def\boA{{\mathcal A}}

\def\boC{{\mathcal C}}
\def\boD{{\mathcal D}}

\def\boF{{\mathcal F}}
\def\boG{{\mathcal G}}

\def\boM{{\mathcal M}}
\def\boP{{\mathcal P}}

\def\boW{{\mathcal W}}
\def\cqfd{\hfill$\Box$}
\def\Res{{\,\rm Res}}
\def\Re{{\rm Re}}
\def\Im{{\rm Im}}

\def\wtPhi{\widetilde{\Phi}}
\def\wtf{\widetilde{f}}
\def\wtF{\widetilde{F}}

\def\wtxi{\widetilde{\xi}}
\def\wtM{\widetilde{M}}
\def\wtN{\widetilde{N}}
\def\wtOmega{\widetilde{\Omega}}
\def\wtSigma{\widetilde{\Sigma}}
\def\wtgamma{\widetilde{\gamma}}
\def\wtz{\widetilde{z}}

\def\whPhi{\widehat{\Phi}}
\def\whxi{\widehat{\xi}}
\def\whf{\widehat{f}}

\def\whpsi{\widehat{\psi}}
\def\whDelta{\widehat{\Delta}}
\def\ii{{\rm i}}
\def\x{{\bf x}}
\def\Uni{\mbox{\rm Uni}}
\def\Pos{\mbox{\rm Pos}}
\def\duni{\mbox{\rm uni}}
\def\dpos{\mbox{\rm pos}}
\def\Sym{\mbox{\rm Sym}}
\def\Nor{\mbox{\rm Nor}}

\def\sl{\mathfrak{sl}}
\def\su{\mathfrak{su}}
\def\boW{{\mathcal W}_{\rho}}
\def\boWp{\boW^{\geq 0}}
\def\boWpp{\boW^{>0}}
\def\boWn{\boW^{\leq 0}}
\def\boWnn{\boW^{<0}}
\def\boWo{{\mathcal W}^0}

\renewcommand{\matrix}[1]{\left(\begin{array}{cc} #1\end{array}\right)}
\newcommand{\smallfrac}[2]{\mbox{$\frac{#1}{#2}$}}
\newcommand{\minimatrix}[1]{\left(\begin{smallmatrix}#1\end{smallmatrix}\right)}
\newtheorem{theorem}{Theorem}
\newtheorem{lemma}{Lemma}
\newtheorem{proposition}{Proposition}
\newtheorem{remark}{Remark}

\newtheorem{claim}{Claim}

\newtheorem{definition}{Definition}

\title{Gluing Delaunay ends to minimal $n$-noids using the DPW method}
\author{Martin Traizet}
\begin{document}
\maketitle

\medskip
{\em Abstract: we construct constant mean curvature surfaces in euclidean space by
gluing $n$ half Delaunay surfaces to a non-degenerate minimal $n$-noid, using the DPW method.}
\section{Introduction}
\label{intro}
In \cite{dorfmeister-pedit-wu}, Dorfmeister, Pedit and Wu have shown that surfaces with non-zero constant
mean curvature  (CMC for short) in euclidean space admit a Weierstrass-type
representation, which means that they can be represented in terms of
holomorphic data. This representation is now called the DPW method.
In \cite{nnoids}, we used the DPW method to construct CMC $n$-noids: genus zero, CMC surfaces with $n$ ends Delaunay type ends.
These $n$-noids can be described as a unit sphere with $n$ half Delaunay surfaces
with small necksizes attached at prescribed points. They had already been constructed
by Kapouleas in \cite{kapouleas} using PDE methods.
\medskip

\begin{figure}
\begin{center}
\includegraphics[width=6cm]{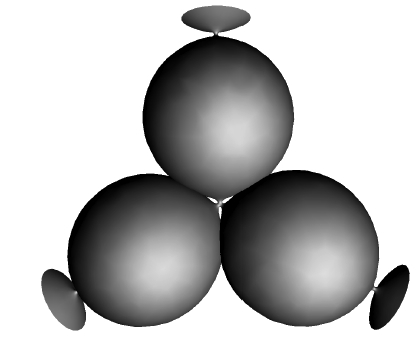}
\hspace{1cm}
\includegraphics[width=5cm]{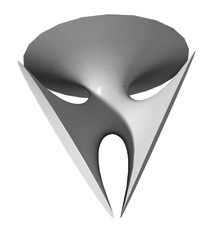}
\end{center}
\caption{A CMC 3-noid (left, image by N. Schmitt \cite{schmitt}) and
a minimal 3-noid (right). There is a tiny ``copy'' of the minimal 3-noid at the center of the CMC 3-noid.}
\label{intro-figure1}
\end{figure}

In the case $n=3$, Alexandrov-embedded CMC trinoids have been classified by
Gro\ss e Brauckman, Kusner and Sullivan in \cite{karsten-kusner-sullivan}.
In particular, equilateral CMC trinoids form a 1-parameter family, parametrized on an open interval.
On one end, equilateral trinoids degenerate like the examples described
above: they look like a sphere with 3 half Delaunay surfaces with small necksizes
attached at the vertices of a spherical equilateral triangle.
On the other end, equilateral trinoids limit, after suitable blow-up, to a minimal
$3$-noid: a genus zero minimal surface with 3 catenoidal ends
(see Figure \ref{intro-figure1}).
\medskip

It seems natural to ask if one can generalize this observation and construct CMC
 $n$-noids by gluing half Delaunay
surfaces with small necksizes to a minimal $n$-noid.
This is indeed the case, and has been done by Mazzeo and Pacard in
\cite{mazzeo-pacard} using PDE methods.
In this paper, we propose a quite simple and natural DPW potential to
construct these examples. We prove:
\begin{theorem}
\label{intro-theorem1}
Let $n\geq 3$ and
let $M_0$ be a non-degenerate minimal $n$-noid.
There exists a smooth family of CMC surfaces $(M_t)_{0<|t|<\epsilon}$
with the following properties:
\begin{enumerate}
\item $M_t$ has genus zero and $n$ Delaunay ends.
\item $\frac{1}{t}M_t$ converges to $M_0$ as $t\to 0$.
\item If $M_0$ is Alexandrov-embedded, all ends of $M_t$ are of unduloid
type if $t>0$ and of nodoid type if $t<0$. Moreover, $M_t$ is Alexandrov-embedded if $t>0$.
\end{enumerate}
\end{theorem}
Non-degeneracy of a minimal $n$-noid will be defined in Section \ref{minimal}.
The two surfaces $M_t$ and $M_{-t}$ are geometrically different: if $M_t$ has an end of unduloid type then the corresponding end of $M_{-t}$ is of nodoid type.
See Proposition \ref{ends-prop1} for more details.
Of course, a minimal $n$-noid is never embedded if $n\geq 3$
so the surfaces $M_t$ are not embedded.
Alexandrov-embedded minimal $n$-noids whose ends have coplanar axes have been classified by Cosin and Ros in \cite{cosin-ros}, and Alexandrov-embedded
CMC $n$-noids whose ends have coplanar axes have been classified by
Gro\ss e-Brauckmann, Kusner and Sullivan in \cite{karsten-kusner-sullivan2}.
\medskip

As already said, these surfaces have already been constructed in \cite{mazzeo-pacard}.
Our motivation to construct them with the DPW method is to answer the following
questions:
\begin{enumerate}
\item How can we produce a DPW potential from the Weierstrass data $(g,\omega)$
of the minimal $n$-noid $M_0$ ?
\item How can we prove, with the DPW method, that $\frac{1}{t}M_t$ converges
to $M_0$ ?
\end{enumerate}
The answer to Question 2 is Theorem \ref{theorem-blowup} in Section \ref{section-blowup},
a general blow-up result in the context of the DPW method.
In \cite{traizet3}, we use the DPW method to construct higher genus CMC surfaces with small necks. Theorem \ref{theorem-blowup} is used to ensure that the necks have asymptotically catenoidal shape.
\section{Non-degenerate minimal $n$-noids}
\label{minimal}
A minimal $n$-noid is a complete, immersed minimal surface in $\R^3$ with
genus zero and $n$ catenoidal ends.
Let $M_0$ be a minimal $n$-noid and $(\Sigma,g,\omega)$ its Weierstrass data.
This means that $M_0$ is parametrized on $\Sigma$ by the Weierstrass Representation
formula:
\begin{equation}
\label{minimal-eq1}
\psi(z)=\Re\int_{z_0}^z\left(\smallfrac{1}{2}(1-g^2)\omega,\smallfrac{\ii}{2}(1+g^2)\omega,g\omega\right)
\end{equation}
Without loss of generality, we can assume that
$\Sigma=\C\cup\{\infty\}\setminus\{p_1,\cdots,p_n\}$, where $p_1,\cdots,p_n$
are complex numbers and $g\neq 0,\infty$ at $p_1,\cdots,p_n$ (by rotating $M_0$ if necessary).
Then $\omega$ needs a double pole at $p_1,\cdots,p_n$ so has $2n-2$ zeros, counting multiplicity.
Since $\omega$ needs a zero at each pole of $g$, with twice the multiplicity, it follows
that $g$ has $n-1$ poles so has degree $n-1$. Hence we may write
\begin{equation}
\label{minimal-eq2}
g=\frac{A(z)}{B(z)}\quad\mbox{ and }\quad
\omega=\frac{B(z)^2\,dz}{\prod_{i=1}^n(z-p_i)^2}
\end{equation}
where
$$
A(z)=\sum_{i=1}^{n}a_i z^{n-i}
\quad\mbox{ and }\quad
 B(z)=\sum_{i=1}^{n}b_i z^{n-i}.
$$
We are going to deform this Weierstrass data, so we see $a_i$, $b_i$ and $p_i$ for
$1\leq i\leq n$ as complex parameters.
We denote by $\x\in\C^{3n}$ the vector of these parameters, and by
$\x_0$ the value of the parameters corresponding to the minimal $n$-noid $M_0$.
\medskip

Let $\gamma_i$ be the homology class of a small circle centered at $p_i$ and define the following
periods for $1\leq i\leq n$ and $0\leq k\leq 2$, depending on the parameter vector $\x\in\C^{3n}$:
$$P_{i,k}(\x)=\int_{\gamma_i}g^k\omega$$
$$P_i(\x)=(P_{i,0}(\x),P_{i,1}(\x),P_{i,2}(\x))\in\C^3$$
$$Q_i(\x)=\int_{\gamma_i}\left(\smallfrac{1}{2}(1-g^2)\omega,\smallfrac{\ii}{2}(1+g^2)\omega,g\omega\right)\in\C^3.$$
Then
$$Q_i(\x)=\left(\smallfrac{1}{2}(P_{i,0(\x)}-P_{i,2}(\x)),\smallfrac{\ii}{2}(P_{i,0}(\x)+P_{i,2}(\x)),P_{i,1}(\x)\right).$$
The components of $Q_i(\x_0)$ are imaginary
because the Period Problem is solved for $M_0$. This gives
\begin{equation}
\label{minimal-eq4}
P_{i,2}(\x_0)=\overline{P_{i,0}(\x_0)}
\quad\mbox{ and } \quad
P_{i,1}(\x_0)\in\ii\R
\end{equation}
Moreover, $\Im(Q_i(\x_0))=-\phi_i$ where $\phi_i$ is the flux vector of $M_0$ at the end $p_i$.
By the Residue Theorem, we have for all $\x$ in a neighborhood of $\x_0$:
$$\sum_{i=1}^n P_i(\x)=0$$
Let $P=(P_1,\cdots,P_{n-1})$ and $Q=(Q_1,\cdots,Q_{n-1})$.
\begin{definition}
\label{minimal-def1}
$M_0$ is non-degenerate if the differential of
$P$ (or equivalently, $Q$) at $\x_0$ has (complex) rank $3n-3$.
\end{definition}
\begin{remark}
\label{minimal-remark1}
If $n\geq 3$, we may (using M\"obius transformations of the sphere) fix the
value of three points, say $p_1,p_2,p_3$. Then "non-degenerate" means that
the differential of $P$ with respect to the remaining parameters is an
isomorphism of $\C^{3n-3}$.
\end{remark}
This notion is related to another standard notion of non-degeneracy:
\begin{definition}
\label{minimal-def2}
$M_0$ is non-degenerate if its space of bounded Jacobi fields has (real) dimension 3.
\end{definition}
\begin{theorem}
\label{minimal-theorem1}
If $M_0$ is non-degenerate in the sense of Definition \ref{minimal-def2},
then $M_0$ is non-degenerate in the sense of Definition \ref{minimal-def1}.
\end{theorem}
Proof. Assume $M_0$ is non-degenerate in the sense of Definition \ref{minimal-def2}.
Then in a neighborhood of $M_0$,
the space $\boM$ of minimal $n$-noids (up to translation) is a smooth manifold of dimension $3n-3$ by
a standard application of the Implicit Function Theorem.
Moreover, if we write $\phi_i\in\R^3$ for the flux vector at the $i$-th end,
then the map $\phi=(\phi_1,\cdots,\phi_n)$ provides a local diffeomorphism between
$\boM$
and the space $V$ of vectors $v=(v_1,\cdots,v_n)\in(\R^3)^n$
such that $\sum_{i=1}^n v_i=0$.
(All this is proved in Section 4 of \cite{cosin-ros} in the case where all ends are coplanar. The argument goes through in the general case.)
Hence given a vector $v\in V$, there exists a deformation $M_t$ of
$M_0$ such that $M_t\in\boM$ and $\frac{d}{dt}\phi(M_t)|_{t=0}=v$.
We may write the Weierstrass data of $M_t$ as above and obtain a
set of parameters $\x(t)$, depending smoothly on $t$, such that
$\x(0)=\x_0$.
Then $dQ(\x_0)\cdot \x'(0)=-\ii v$.
Since $Q$ is holomorphic, its differential is complex-linear so $dQ(\x_0)$
has complex rank equal to $\dim V=3n-3$.\cqfd

\medskip

If all ends of $M_0$ have coplanar axes, then $M_0$ is
non-degenerate in the sense of Definition \ref{minimal-def2}
by Proposition 2 in  \cite{cosin-ros}.
In particular, the (most symmetric) $n$-noids of Jorge-Meeks are non-degenerate.
This implies that generic $n$-noids in the component of the Jorge-Meeks $n$-noid
are non-degenerate.
\section{Background}
\label{background}
In this section, we recall standard notations and results used in the
DPW method. We work in the ``untwisted'' setting.
\subsection{Loop groups}
\label{background-section1}
A loop is a smooth map from the unit circle $\S^1=\{\lambda\in\C:|\lambda|=1\}$ to a matrix group.
The circle variable is denoted $\lambda$ and called the spectral parameter.
For $\rho>0$, we denote $\D_{\rho}=\{\lambda\in\C\,:\, |\lambda|<\rho\}$,
$\D_{\rho}^*=\D_{\rho}\setminus\{0\}$ and $\D=\D_1$.
\begin{itemize}
\item If $G$ is a matrix Lie group (or Lie algebra), $\Lambda G$ denotes the group
(or algebra)
of smooth maps $\Phi:\S^1\to G$.
\item $\Lambda_+ SL(2,\C)\subset\Lambda SL(2,\C)$ is the subgroup of maps $B$ which extend holomorphically to $\D$ with $B(0)$ upper triangular.
\item $\Lambda_+^{\R}SL(2,\C)\subset\Lambda_+ SL(2,\C)$ is the subgroup of maps $B$ such that
$B(0)$ has positive entries on the diagonal.
\end{itemize}
\begin{theorem}[Iwasawa decomposition]
The multiplication $\Lambda SU(2)\times \Lambda_+^{\R} SL(2,\C)\to\Lambda SL(2,\C)$
is a diffeomorphism. The unique splitting of an element $\Phi\in\Lambda SL(2,\C)$ as
$\Phi=FB$ with $F\in\Lambda SU(2)$ and $B\in\Lambda_+^{\R} SL(2,\C)$
is called Iwasawa decomposition.
$F$ is called the unitary factor of $\Phi$ and
denoted $\Uni(\Phi)$.
$B$ is called the positive factor and denoted $\Pos(\Phi)$.
\end{theorem}
\subsection{The matrix model of $\R^3$}
\label{background-section2}
In the DPW method, one identifies $\R^3$ 
with the Lie algebra $\su(2)$ by
$$x=(x_1,x_2,x_3)\in\R^3\longleftrightarrow X=-\ii\matrix{
-x_3&x_1+\ii x_2\\x_1-\ii x_2 &x_3}\in\mathfrak{su}(2).$$
We have $\det(X)=\|x\|^2$.
The group $SU(2)$ acts as linear isometries on $\su(2)$ by
conjugation: $H\cdot X=HXH^{-1}$.
\subsection{The DPW method}
\label{background-section3}
The input data for the DPW method is a quadruple $(\Sigma,\xi,z_0,\phi_0)$ where:
\begin{itemize}
\item $\Sigma$ is a Riemann surface.
\item
$\xi=\xi(z,\lambda)$ is a $\Lambda\sl(2,\C)$-valued holomorphic 1-form on $\Sigma$
called the DPW potential. More precisely,
\begin{equation}
\label{equation-xi}
\xi=\matrix{\alpha&\lambda^{-1}\beta\\ \gamma&-\alpha}
\end{equation}
where $\alpha(z,\lambda)$, $\beta(z,\lambda)$, $\gamma(z,\lambda)$
are holomorphic 1-forms on $\Sigma$ with respect to the $z$ variable,
and are holomorphic with respect to $\lambda$ in the disk $\D_{\rho}$
for some $\rho>1$.
\item $z_0\in\Sigma$ is a base point.
\item $\phi_0\in\Lambda SL(2,\C)$ is an initial condition.
\end{itemize}
Given this data, the DPW method is the following procedure.
\begin{itemize}
\item Let $\wtSigma$ be the universal cover of $\Sigma$ and $\wtz_0\in\wtSigma$
be an arbitrary element in the fiber of $z_0$. Solve the Cauchy Problem on $\wtSigma$:
\begin{equation}
\label{cauchy-problem}
\left\{\begin{array}{l}
d\Phi(z,\lambda)=\Phi(z,\lambda)\xi(z,\lambda)\\
\Phi(\wtz_0,\lambda)=\phi_0(\lambda)\end{array}\right.\end{equation}
to obtain a solution
$\Phi: \wtSigma\to\Lambda SL(2,\C)$.
\item Compute the Iwasawa decomposition $(F(z,\cdot),B(z,\cdot))$ of $\Phi(z,\cdot)$.
\item Define $f:\wtSigma\to\su(2)\sim\R^3$ by the Sym-Bobenko formula:
\begin{equation}
\label{sym-bobenko}
f(z)=-2\ii\frac{\partial F}{\partial\lambda}(z,1)F(z,1)^{-1}=:\Sym(F(z,\cdot)).
\end{equation}
Then $f$ is a CMC-1 (branched) conformal immersion.
$f$ is regular at $z$ (meaning unbranched) if and only if $\beta(z,0)\neq 0$.
Its Gauss map is given by
\begin{equation}
\label{eq-N}
N(z)=-\ii F(z,1)\matrix{-1 &0 \\ 0 & 1}F(z,1)^{-1}=:\Nor(F(z,\cdot)).
\end{equation}
The differential of $f$ satisfies
\begin{equation}
\label{eq-df}
df(z)=2\ii\,B_{11}(z,0)^2F(z,1)\matrix{0&\beta(z,0)\\\overline{\beta}(z,0)&0}F(z,1)^{-1}.
\end{equation}
\end{itemize}
\begin{remark}
\label{remark-signs}
\begin{enumerate}
\item In \cite{nnoids}, I have opposite signs in Equations \eqref{sym-bobenko} and
\eqref{eq-N}. This is unfortunate because it makes the basis $(f_x,f_y,N)$ negatively oriented.
Equation \eqref{sym-bobenko} is the right formula, which one obtains by untwisting the standard Sym-Bobenko
formula in the twisted case. See \cite{kilian-rossman-schmitt} or \cite{schmitt-kilian-kobayashi-rossman}.
\item I have not been able to find Formula \eqref{eq-df} in the litterature. Of course, the DPW method constructs a moving frame for $f$, so one has a formula for $df$, but usually it is written in a special coordinate system and only in the ``twisted'' setting. For the interested reader, I derive Equation \eqref{eq-df} from the Sym-Bobenko formula at the end of Appendix \ref{appendix-iwasawa}.
\end{enumerate}
\end{remark}
\subsection{The Monodromy Problem}
Assume that $\Sigma$ is not simply connected so its universal cover $\wtSigma$
is not trivial.
Let $\mbox{Deck}(\wtSigma/\Sigma)$ be the group of fiber-preserving diffeomorphisms of $\wtSigma$.
For $\gamma\in\mbox{Deck}(\wtSigma/\Sigma)$, let
$$\boM_{\gamma}(\Phi)(\lambda)=\Phi(\gamma(z),\lambda)\Phi(z,\lambda)^{-1}$$ 
be the monodromy of $\Phi$ with respect to $\gamma$
(which is independent of $z\in\wtSigma$).
The standard condition which ensures that the immersion $f$ descends to
a well defined immersion on $\Sigma$
is the following system of equations, called the Monodromy Problem.
\begin{equation}
\label{monodromy-problem}
\forall \gamma\in\mbox{Deck}(\wtSigma/\Sigma)\quad
\left\{\begin{array}{lc}
\boM_{\gamma}(\Phi)\in\Lambda SU(2)\qquad &(i)\\
\boM_{\gamma}(\Phi)(1)=\pm I_2\qquad&(ii)\\
\frac{\partial\boM_{\gamma}(\Phi)}{\partial\lambda}(1)=0\qquad &(iii)
\end{array}\right.\end{equation}
One can identify $\mbox{Deck}(\wtSigma/\Sigma)$ with the fundamental group
$\pi_1(\Sigma,z_0)$ (see for example Theorem 5.6 in \cite{forster}), so we will in general see $\gamma$ as an element of
$\pi_1(\Sigma,z_0)$. Under this identification, the monodromy of $\Phi$ with
respect to $\gamma\in\pi_1(\Sigma,z_0)$ is given by
$$\boM_{\gamma}(\Phi)(\lambda)=\Phi(\wtgamma(1),\lambda)\Phi(\wtgamma(0),\lambda)^{-1}$$
where $\wtgamma:[0,1]\to\wtSigma$ is the lift of $\gamma$ such that $\wtgamma(0)=\wtz_0$.
\subsection{Gauging}
\label{background-section6}
\begin{definition}
A gauge on $\Sigma$ is a map $G:\Sigma\to\Lambda_+ SL(2,\C)$ such that
$G(z,\lambda)$ depends holomorphically on $z\in\Sigma$ and
$\lambda\in\D_{\rho}$ for some $\rho>1$.
\end{definition}
Let $\Phi$ be a solution of $d\Phi=\Phi\xi$ and $G$ be a gauge.
Let $\whPhi=\Phi\times G$. Then $\whPhi$ and $\Phi$
define the same immersion $f$. This is called ``gauging''.
The gauged potential is
$$\whxi=\whPhi^{-1}d\whPhi=G^{-1}\xi G+G^{-1}dG$$
and will be denoted $\xi\cdot G$, the dot denoting the action of the gauge
group on the potential.
\subsection{Functional spaces}
\label{section-functional-spaces}
We decompose a smooth function $f:\S^1\to\C$ in Fourier
series
$$f(\lambda)=\sum_{i\in\Z}f_i\lambda^i$$
Fix some $\rho>1$ and define
$$\|f\|=\sum_{i\in\Z}|f_i|\rho^{|i|}$$
Let $\boW$ be the space of functions $f$ with finite norm.
This is a Banach algebra, classically called the Wiener algebra when $\rho=1$
(see Proposition \ref{appendix-prop-algebra} in appendix \ref{appendix-iwasawa}).
Functions in $\boW$ extend holomorphically to the annulus $\frac{1}{\rho}<|\lambda|<\rho$.
\medskip

We define $\boWp$, $\boWpp$, $\boWn$ and $\boWnn$ as the subspaces of functions $f$ such that $f_i=0$
for $i<0$, $i\leq 0$, $i>0$ and $i\geq 0$, respectively.
Functions in $\boWp$ extend holomorphically to the disk $\D_{\rho}$ and
satisfy $|f(\lambda)|\leq \|f\|$ for all $\lambda\in\D_{\rho}$.
We write $\boWo\sim\C$ for the subspace of constant functions, so we have
a direct sum $\boW=\boWnn\oplus\boWo\oplus\boWpp$. A function $f$
will be decomposed as $f=f^-+f^0+f^+$ with
$(f^-,f^0,f^+)\in\boWnn\times\boWo\times\boWpp$.
\medskip

We define the star operator by
$$f^*(\lambda)=\overline{f\left(\frac{1}{\overline{\lambda}}\right)}=\sum_{i\in\Z}\overline{f_{-i}}\lambda^i$$
The involution $f\mapsto f^*$ exchanges $\boWp$ and $\boWn$.
We have $\lambda^*=\lambda^{-1}$ and $c^*=\overline{c}$ if $c$ is a constant.
A function $f$ is real on the unit circle if and only if $f=f^*$.
\medskip

If $L$ is a loop group, we denote $L_{\rho}\subset L$ the subgroup
of loops whose entries are in $\boW$.
If $\Phi\in\Lambda SL(2,\C)_{\rho}$ and $(F,B)$ is its Iwasawa decomposition, then in
fact $F\in\Lambda SU(2)_{\rho}$ and $B\in\Lambda_+^{\R} SL(2,\C)_{\rho}$
(see Proposition \ref{appendix-iwasawa-proposition1} in Appendix \ref{appendix-iwasawa}).
Moreover, Iwasawa decomposition is smooth, as a map between Banach manifolds
(see Theorem \ref{appendix-iwasawa-theorem} in Appendix \ref{appendix-iwasawa}).
\section{A blow-up result}
In this section, we consider a one-parameter family of DPW potential $\xi_t$ with solution
$\Phi_t$ and assume that $\Phi_0(z,\lambda)$ is independent of $\lambda$. Then
its unitary part $F_0(z,\lambda)$ is independent of $\lambda$. The Sym Bobenko formula yields
that $f_0\equiv 0$, so the family $f_t$ collapses to the origin as $t=0$. The following theorem
says that the blow-up $\frac{1}{t}f_t$ converges to a minimal surface whose Weierstrass data is
explicitly computed.
\label{section-blowup}
\begin{theorem}
\label{theorem-blowup}
Let $\Sigma$ be a Riemann surface, $(\xi_t)_{t\in I}$ a family of DPW potentials on
$\Sigma$ and $(\Phi_t)_{t\in I}$ a family of solutions of $d\Phi_t=\Phi_t\xi_t$ on the universal cover $\wtSigma$ of $\Sigma$, where $I\subset\R$ is a neighborhood of $0$.
Fix a base point $z_0\in\wtSigma$.
Assume that
\begin{enumerate}
\item $(t,z)\mapsto \xi_t(z,\cdot)$ and $t\mapsto \Phi_t(z_0,\cdot)$ are $C^1$ maps into
$\Lambda\sl(2,\C)_{\rho}$ and $\Lambda SL(2,\C)_{\rho}$, respectively.
\item For all $t\in I$, $\Phi_t$ solves the Monodromy Problem \eqref{monodromy-problem}.
\item $\Phi_0(z,\lambda)$ is independent of $\lambda$:
$$\Phi_0(z,\lambda)=\matrix{a(z)&b(z)\\c(z)&d(z)}$$
\end{enumerate}
Let $f_t=\Sym(\Uni(\Phi_t)):\Sigma\to\R^3$ be the CMC-1 immersion given
by the DPW method.
Then 
$$\lim_{t\to 0}\frac{1}{t}f_t(z)=\psi(z)$$ where
$\psi:\Sigma\to\R^3$ is a (possibly branched) minimal immersion with the following Weierstrass data:
$$g(z)=\frac{-a(z)}{c(z)}\quad\mbox{ and }\quad
\omega=4c(z)^2\frac{\partial \xi_{t;12}^{(-1)}}{\partial t}|_{t=0}.$$
The limit is for the uniform $C^1$ convergence on compact subsets of $\Sigma$.
\end{theorem}
Here $\xi_{t;12}^{(-1)}$ denotes the coefficient of $\lambda^{-1}$ in the upper right entry of $\xi_t$. In case $\omega=0$, the minimal immersion degenerates into a point and $\psi$ is constant.
\medskip

Proof: by standard ODE theory, $(t,z)\mapsto \Phi_t(z,\cdot)$ is a $C^1$ map
into $\Lambda SL(2,\C)_{\rho}$.
Let $(F_t,B_t)$ be the Iwasawa decomposition of $\Phi_t$.
By Theorem \ref{appendix-iwasawa-theorem} in Appendix \ref{appendix-iwasawa}, $(t,z)\mapsto F_t(z,\cdot)$ and
$(t,z)\mapsto B_t(z,\cdot)$ are $C^1$ maps into $\Lambda SU(2)_{\rho}$
and $\Lambda_+^{\R} SL(2,\C)_{\rho}$, respectively.
At $t=0$, $\Phi_0$ is constant with respect to $\lambda$, so its Iwasawa decomposition is the standard $QR$ decomposition:
$$F_0=\frac{1}{\sqrt{|a|^2+|c|^2}}\matrix{a&-\overline{c}\\c&\overline{a}}\qquad
B_0=\frac{1}{\sqrt{|a|^2+|c|^2}}\matrix{|a|^2+|c|^2&\overline{a}b+\overline{c}d\\0&1}.$$
The Sym-Bobenko formula \eqref{sym-bobenko} yields $f_0=0$.
Let $\mu_t=B_{t;11}^0$ and $\beta_t=\xi_{t;12}^{(-1)}$.
By Equation \eqref{eq-df}, we have
$$df_t(z)=2\ii\,\mu_t(z)^2F_t(z,1)\matrix{0&\beta_t(z)\\\overline{\beta}_t(z)&0}F_t(z,1)^{-1}.$$
Hence $(t,z)\mapsto df_t(z)$ is a $C^1$ map.
At $t=0$, $\xi_0$ is constant with respect to $\lambda$, so $\beta_0=0$.
Define $\wtf_t(z)=\frac{1}{t}f_t(z)$ for $t\neq 0$. Then $d\wtf_t(z)$ extends at $t=0$,
as a continous function of $(t,z)$, by
\begin{eqnarray*}
d\wtf_0=\frac{d}{dt}df_t|_{t=0}
&=&2\ii \matrix{a&-\overline{c}\\c&\overline{a}}\matrix{0&\beta'\\\overline{\beta'}&0}\matrix{\overline{a}&\overline{c}\\-c&a}\\
&=&2\ii\matrix{-ac\beta'-\overline{ac\beta'}&a^2\beta'-\overline{c^2\beta'}\\\overline{a^2\beta'}-c^2\beta'&ac\beta'+\overline{ac\beta'}}
\end{eqnarray*}
where $\beta'=\frac{d}{dt}\beta_t|_{t=0}$.
In euclidean coordinates, this gives
$$d\wtf_0=4\;\Re\left[\smallfrac{1}{2}(c^2-a^2)\beta',\smallfrac{\ii}{2}(c^2+a^2)\beta',-ac\beta'\right].$$
Writing $g=\frac{-a}{c}$ and
$\omega=4c^2\beta'$, we obtain
$$\wtf_0(z)=\wtf_0(z_0)+\Re\int_{z_0}^z\left[\smallfrac{1}{2}(1-g^2)\omega,
\smallfrac{\ii}{2}(1+g^2)\omega,g\omega\right]$$
and we see that $\wtf_0$ is a minimal surface with Weierstrass data
$(g,\omega)$.
The last statement of Theorem \ref{theorem-blowup} comes from the fact that $d\wtf_t$ converges uniformly to $d\wtf_0$ on compact subsets of $\Sigma$.

\subsection{Example}
As an example, we consider the family of Delaunay surfaces given by the following DPW potential
in $\C^*$:
$$\xi_t(z,\lambda)=\matrix{0&\lambda^{-1}r+s\\\lambda r+s&0}\frac{dz}{z}
\quad\mbox{ with }\quad
\left\{\begin{array}{l}
r+s=\frac{1}{2}\\
rs=t\\
r<s\end{array}\right.$$
with initial condition $\Phi_t(1)=I_2$.
As $t\to 0$, we have $(r,s)\to(0,\frac{1}{2})$. We have
$$\Phi_0(z,\lambda)=\exp\matrix{0&\frac{1}{2}\\\frac{1}{2}&0}\log z=\frac{1}{2\sqrt{z}}\matrix{z+1&z-1\\z-1&z+1}$$
$$\frac{\partial\xi_t}{\partial t}|_{t=0}=\matrix{0&2\lambda^{-1}\\2\lambda&0}\frac{dz}{z}.$$
Theorem \ref{theorem-blowup} applies and gives
$$g(z)=\frac{1+z}{1-z}\quad\mbox{ and }\quad
\omega(z)=4\left(\frac{z-1}{2\sqrt{z}}\right)^2\frac{2\,dz}{z}=2\left(\frac{z-1}{z}\right)^2dz.$$
This is the Weierstrass data of a horizontal catenoid of waist-radius 4 and axis $Ox_1$, with $x_1\to +\infty$ at the end $z=0$.
\section{The DPW potential}
We now start the proof of Theorem \ref{intro-theorem1}.
Let $(g,\omega)$ be the Weierstrass data of the given minimal $n$-noid $M_0$,
written as in Section \ref{minimal}.
We introduce $3n$ $\lambda$-dependent parameters $a_i$, $b_i$ and $p_i$ for $1\leq i\leq n$ in the functional space $\boWp$.
The vector of these parameters is denoted $\x\in(\boWp)^{3n}$.
The parameter $\x$ is in a neighborhood of a (constant) central value $\x_0\in(\boWo)^{3n}$ which correspond to the Weierstrass data of $M_0$, written as in Section \ref{minimal}.
We define
$$A_{\x}(z,\lambda)=\sum_{i=1}^n a_i(\lambda) z^{n-i}$$
$$B_{\x}(z,\lambda)=\sum_{i=1}^n b_i(\lambda) z^{n-i}$$
\begin{equation}
\label{eq-g}
g_{\x}(z,\lambda)=\frac{A_{\x}(z,\lambda)}{B_{\x}(z,\lambda)}
\end{equation}
\begin{equation}
\label{eq-omega}
\omega_{\x}(z,\lambda)=\frac{B_{\x}(z,\lambda)^2\,dz}{\prod_{i=1}^n(z-p_i(\lambda))^2}.
\end{equation}
For $t$ in a neighborhood of $0$ in $\R$, we consider the following DPW potential:
$$\xi_{t,\x}(z,\lambda)=\matrix{0 & \smallfrac{1}{4}t(\lambda-1)^2\lambda^{-1}\omega_{\x}(z,\lambda)\\ dg_{\x}(z,\lambda) & 0}.$$
We fix a base point $z_0$, away from the poles of $g$ and $\omega$, and we take the initial condition
$$\phi_0(\lambda)=\matrix{g_{\x}(z_0,\lambda)&1\\-1&0}.$$
These choices are motivated by the following observations:
\begin{enumerate}
\item At $t=0$, we have
$$\xi_{0,\x}(z,\lambda)=\matrix{0&0\\dg_{\x}(z,\lambda)&0}.$$
The solution of the Cauchy Problem \eqref{cauchy-problem} is given by
\begin{equation}
\label{eq-Phi0}
\Phi_{0,\x}(z,\lambda)=\matrix{g_{\x}(z,\lambda)&1\\-1&0}
\end{equation}
which is well-defined, so the Monodromy Problem \eqref{monodromy-problem} is solved at $t=0$.
\item The same conclusion holds if $\lambda=1$ instead of $t=0$. In particular,
Items (ii) and (iii) of the Monodromy Problem \eqref{monodromy-problem} are automatically solved.
\item At $\x=\x_0$, we have $g_{\x_0}=g$ so $\Phi_{0,\x_0}(z,\lambda)$ is independent of $\lambda$.
Moreover,
$$\frac{\partial \xi_{t,\x_0;12}^{(-1)}}{\partial t}|_{t=0}=\frac{\omega}{4}.$$
Provided the Monodromy Problem is solved for all $t$ in a neighborhood of $0$, Theorem \ref{theorem-blowup} applies and the limit minimal surface has Weierstrass data
$(g,\omega)$ so is the minimal $n$-noid $M_0$, up to translation (see details in Section \ref{section-convergence}).
\end{enumerate}
\subsection{Regularity}
\label{section-regularity}
Our potential $\xi_{t,\x}$ has poles at the zeros of $B_{\x}$ and the points
$p_1,\cdots,p_n$. (At $\infty$, we have $\omega_{\x}\sim b_1^2z^{-2}dz$ which is holomorphic.)
We want the zeros of $B_{\x}$ to be apparent singularities, so we require the
potential to be gauge-equivalent to a regular potential in a neighborhood of these
points.
Consider the gauge
$$G_{\x}(z,\lambda)=\matrix{g_{\x}(z,\lambda)^{-1} & -1\\ 0 & g_{\x}(z,\lambda)}$$
The gauged potential is
$$\whxi_{t,\x}:=\xi_{t,\x}\cdot G_{\x}=\matrix{0 & \smallfrac{1}{4}t(\lambda-1)^2\lambda^{-1} g_{\x}^2 \omega_{\x} \\ g_{\x}^{-2}dg_{\x} & 0}.$$
We have
$$g_{\x}^{-2}dg_{\x}=\frac{A_{\x}'B_{\x}-A_{\x}B_{\x}'}{A_{\x}^2}\quad\mbox{ and }\quad
g_{\x}^2\omega_{\x}=\frac{A_{\x}^2\,dz}{\prod_{i=1}^n(z-p_i)^2}.$$
Let $\zeta$ be a zero of $B_{\x_0}$ (recall that $B_{\x_0}$ does not depend on $\lambda$).
Then $A_{\x_0}(\zeta)\neq 0$.
By continuity, there exists a neighborhood $U$ of $\zeta$ such that for $z\in U$, $\lambda\in\D_{\rho}$ and $\x$ close enough to $\x_0$, $A_{\x}(z,\lambda)\neq 0$.
So $\whxi_{t,\x}$ is holomorphic in $U\times\D_{\rho}^*$ and moreover, $\whxi_{t,\x;12}^{(-1)}\neq 0$.
This ensures that the immersion extends analytically to $U$ and is unbranched in $U$.
\section{The monodromy problem}
\label{monodromy}
\subsection{Formulation of the problem}
\label{monodromy-section1}
For $i\in[1,n]$, we denote $p_{i,0}$ the central value of the parameter $p_i$
(so $p_{1,0},\cdots,p_{n,0}$ are the ends of the minimal $n$-noid $M_0$).
We consider the following $\lambda$-independent domain on the Riemann sphere:
\begin{equation}
\label{eq-Omega}
\Omega=\{z\in\C:\forall i\in[1,n],|z-p_{i,0}|>\varepsilon\}\cup\{\infty\}
\end{equation}
where $\varepsilon>0$ is a fixed, small enough number such that the disks $D(p_{i,0},8\varepsilon)$
for $1\leq i\leq n$ are disjoint.
As in \cite{nnoids}, we first construct a family of immersions $f_t$ on $\Omega$. Then we extend
$f_t$ to an $n$-punctured sphere  in Proposition \ref{prop-extends}.
\medskip

Let $\wtOmega$ be the universal cover of $\Omega$ and
$\Phi_{t,\x}(z,\lambda)$ be the solution of the following Cauchy
Problem on $\wtOmega$:
\begin{equation}
\left\{\begin{array}{l}
d\Phi_{t,\x}(z,\lambda)=\Phi_{t,\x}(z,\lambda)\xi_{t,\x}(z,\lambda)\\
\Phi_{t,\x}(\wtz_0,\lambda)=\phi_0\end{array}\right.
\end{equation}
We denote $\gamma_1,\cdots,\gamma_{n-1}$ a set of generators of the fundamental group
$\pi_1(\Omega,z_0)$, with $\gamma_i$ encircling the point $p_{i,0}$.
We may assume that each $\gamma_i$ is represented by a fixed curve avoiding the poles of $\xi_{t,\x}$.
Let
$$M_i(t,\x)=\boM_{\gamma_i}(\Phi_{t,\x})$$
be the monodromy of $\Phi_{t,\x}$ along $\gamma_i$.
By Equation \eqref{eq-Phi0}, we have $M_i(0,\x)=I_2$.
Recall that the matrix exponential is a local diffeomorphism from a neighborhood of $0$
in the Lie algebra $\sl(2,\C)$ (respectively $\su(2)$) to
a neighborhood of $I_2$ in $SL(2,\C)$ (respectively $SU(2)$).
The inverse diffeomorphism is denoted $\log$.
For $t\neq 0$ small enough and $\lambda\in\D_{\rho}\setminus\{1\}$, we define as in \cite{nnoids}
$$\wtM_i(t,\x)(\lambda)=\frac{4\lambda}{t(\lambda-1)^2}\log M_i(t,\x)(\lambda).$$
\begin{proposition}
\label{monodromy-prop1}
\begin{enumerate}
\item $\wtM_i(t,\x)(\lambda)$ extends smoothly at $t=0$ and $\lambda=1$, and 
each entry $\wtM_{i;k\ell}$ is a smooth map from a neighborhood of
$(0,\x_0)$ in $\R\times(\boWp)^3$ to $\boW$.
\item At $t=0$, we have
\begin{equation}
\label{monodromy-eq2}
\wtM_i(0,\x)(\lambda)=\matrix{\boP_{i,1}(\x) & \boP_{i,2}(\x) \\ -\boP_{i,0}(\x) & -\boP_{i,1}(\x)}
\quad\mbox{ where }
\boP_{i,k}(\x)=\int_{\gamma_i}g_{\x}^k\omega_{\x}.
\end{equation}
\item The Monodromy Problem \eqref{monodromy-problem} is equivalent to
\begin{equation}
\label{monodromy-eq3}
\wtM_i(t,\x)\in\Lambda\su(2)\quad\mbox{ for }\quad 1\leq i\leq n-1.
\end{equation}
\end{enumerate}
\end{proposition}
Proof: we follow the proof of Proposition 1 in \cite{nnoids}.
We first consider the case where the parameter $\x=(a_i,b_i,p_i)_{1\leq i\leq n}$ is constant with
respect to $\lambda$, so $\x\in\C^{3n}$.
For $(\mu,\x)$ in a neighborhood of $(0,\x_0)$ in $\C\times\C^{3n}$, we define
$$\whxi_{\mu,\x}(z)=\matrix{0&\mu\,\omega_{\x}(z)\\dg_{\x}(z)&0}$$
where $\omega_{\x}$ and $g_{\x}$ are defined by Equations \eqref{eq-g} and \eqref{eq-omega}, except that
$a_i$, $b_i$, $p_i$ are constant complex numbers.
Let $\whPhi_{\mu,\x}$ be the solution of the Cauchy Problem
$d\whPhi_{\mu,\x}=\whPhi_{\mu,\x}\whxi_{\mu,\x}$ in $\wtOmega$ with initial condition
$\whPhi_{\mu,\x}(\wtz_0)=\phi_0$.
Let $N_i(\mu,\x)=\boM_{\gamma_i}(\whPhi_{\mu,\x})$.
By standard ODE theory, each entry of $N_i$ is a holomorphic function of $(\mu,\x)$.
At $\mu=0$, $\whPhi_{0,\x}$ is given by Equation \eqref{eq-Phi0}, so in particular
$N_i(0,\x)=I_2$. Hence
$$\wtN_i(\mu,\x):=\frac{1}{\mu}\log N_i(\mu,\x)$$
extends holomorphically at $\mu=0$ with value
$\wtN_i(0,\x)=\frac{\partial N_i}{\partial\mu}(0,\x)$.
By Proposition 8 in Appendix A of \cite{nnoids},
$$\frac{\partial N_i}{\partial \mu}(0,\x)=\int_{\gamma_i}\whPhi_{0,\x}\frac{\partial\whxi_{\mu,\x}}{\partial\mu}|_{\mu=0}\whPhi_{0,\x}^{-1}.$$
Hence
\begin{equation}
\label{monodromy-eqN}
\wtN_i(0,\x)=\int_{\gamma_i}\matrix{g_{\x} &1\\-1&0}
\matrix{0 &\omega_{\x} \\ 0 & 0}\matrix{0&-1\\ 1&g_{\x}}
=\int_{\gamma_i}\matrix{g_{\x}\omega_{\x} &g_{\x}^2\omega_{\x}\\-\omega_{\x} & -g_{\x}\omega_{\x}}.
\end{equation}
For $(t,\x)$ in a neighborhood of $(0,\x_0)$ in $\R\times(\boWp)^{3n}$, we have
$$\xi_{t,\x}(z,\lambda)=\whxi_{\mu(t,\lambda),\x(\lambda)}(z)
\quad\mbox{ with }\quad\mu(t,\lambda)=\frac{t(\lambda-1)^2}{4\lambda}.$$
Hence
$$
M_i(t,\x)(\lambda)=N_i(\mu(t,\lambda),\x(\lambda))
\quad\mbox{ and }\quad
\wtM_i(t,\x)(\lambda)=\wtN_i(\mu(t,\lambda),\x(\lambda)).$$
By substitution (see Proposition 9 in Appendix B of \cite{nnoids}), each entry of $\wtM_i$ is
is a smooth map from a neighborhood of
$(0,\x_0)$ in $\R\times(\boWp)^3$ to $\boW$. Moreover, $\wtM_i(0,\x)$
is given by Equation \eqref{monodromy-eqN}.
The fact that $\wtM_i$ extends holomorphically at $\lambda=1$ implies that Points (ii) and
(iii) of Problem \eqref{monodromy-problem} are automatically satisfied.
Since $\lambda^{-1}(\lambda-1)^2\in\R$ for $\lambda\in\S^1$,
Equation (i) of Problem \eqref{monodromy-problem} is equivalent to Equation \eqref{monodromy-eq3}.\cqfd
\subsection{Solution of the monodromy problem}
Without loss of generality, we may (using a M\"obius transformation of the sphere) fix the value of $p_1$, $p_2$ and $p_3$. We still denote $\x\in(\boWp)^{3n-3}$ the vector of the remaining
parameters.
\label{monodromy-section3}
\begin{proposition}
\label{monodromy-prop2}
Assume that the given minimal $n$-noid is non-degenerate.
For $t$ in a neighborhood of $0$, there exists a smooth function $\x(t)\in(\boWp)^{3n-3}$ such that $\wtM_i(t,\x(t),\cdot)\in\Lambda\su(2)$
for $1\leq i \leq n-1$.
Moreover, $\x(0)=\x_0$.
\end{proposition}
Proof: 
recalling the definition of $P_{i,k}$ in Section \ref{minimal} and $\boP_{i,k}$ in
Equation \eqref{monodromy-eq2}, we have
$$\boP_{i,k}(\x)(\lambda)=P_{i,k}(\x(\lambda)).$$
Hence $\boP_{i,k}$ is a smooth map from a neighborhood of $\x_0$ in $(\boWp)^{3n-3}$ to $\boWp$. Moreover, since
$\x_0$ is constant, we have for $X\in(\boWp)^{3n-3}$:
\begin{equation}
\label{eq-dboP}
(d\boP_{i,k}(\x_0)X)(\lambda)=dP_{i,k}(\x_0)X(\lambda).
\end{equation}
Let $\boP_i=(\boP_{i,0},\boP_{i,1},\boP_{i,2})$ and $\boP=(\boP_1,\cdots,\boP_{n-1})$.
By the non-degeneracy hypothesis and Remark \ref{minimal-remark1}, $dP(\x_0)$ is an automorphism of $\C^{3n-3}$, so
$d\boP(\x_0)$ is an automorphism of $(\boWp)^{3n-3}$
and restricts to an automorphism of $(\boWpp)^{3n-3}$.
\medskip

We define the following smooth maps with value in $\boW$
(the star operator is defined in Section \ref{section-functional-spaces})
$$\boF_i(t,\x)=\wtM_{i,11}(t,\x)+\wtM_{i,11}(t,\x)^*$$
$$\boG_i(t,\x)=\wtM_{i,12}(t,\x)+\wtM_{i,21}(t,\x)^*$$
Problem \eqref{monodromy-eq3} is equivalent to $\boF_i=\boG_i=0$.
Actually, by definition, $\boF_i=\boF_i^*$, so Problem \eqref{monodromy-eq3} is equivalent to
$$\boF_i(t,\x)^+=0,\quad \Re(\boF_i(t,\x)^0)=0\quad\mbox{ and } \boG_i(t,\x)=0
\quad\mbox{ for $1\leq i\leq n-1$.}$$
At $t=0$, we have by Equation \eqref{monodromy-eq2}:
$$\boF_i(0,\x)=\boP_{i,1}(\x)+\boP_{i,1}(\x)^*$$
$$\boG_i(0,\x)=\boP_{i,2}(\x)-\boP_{i,0}(\x)^*$$
Equation \eqref{minimal-eq4} tells us precisely that
that at the central value, we have $\boF_i(0,\x_0)=0$ and $\boG_i(0,\x_0)=0$.
We have for $X\in(\boWp)^{3n-3}$:
$$d\boF_i(0,\x_0)X=d\boP_{i,1}(x_0)X+(d\boP_{i,1}(\x_0)X)^*$$
$$d\boG_i(0,\x_0)X=d\boP_{i,2}(x_0)X-(d\boP_{i,0}(\x_0)X)^*$$
Projecting on $\boWpp$ and $\boWnn$ we obtain:
$$(d\boF_i(0,\x_0)X)^+=d\boP_{i,1}(\x_0)X^+$$
$$(d\boG_i(0,\x_0)X)^+=d\boP_{i,2}(\x_0)X^+$$
$$(d\boG_i(0,\x_0)X)^-=-(d\boP_{i,0}(\x_0)X^+)^*$$
$$(d\boG_i(0,\x_0)X)^{-*}=-d\boP_{i,0}(\x_0)X^+.$$
Hence the operator
$$\left[d\boF_i(0,\x_0)^+,d\boG_i(0,\x_0)^+,d\boG_i(0,\x_0)^{-*}\right]_
{1\leq i\leq n-1}$$
only depends on $X^+$ and is an automorphism of $(\boWpp)^{3n-3}$
because $d\boP(\x_0)$ is.
Projecting on $\boWo$ we obtain:
$$(d\boF_i(0,\x_0)X)^0=2\,\Re\left(d\boP_{i,1}(\x_0)X^0\right)$$
$$(d\boG_i(0,\x_0)X)^0=d\boP_{i,2}(\x_0)X^0-\overline{d\boP_{i,0}(\x_0)X^0}.$$
Hence the $\R$-linear operator
$$\left[\Re(d\boF_i(0,\x_0)^0),d\boG_i(0,\x_0)^0\right]_{1\leq i\leq n-1}$$
only depends on $X^0$ and is surjective from $\C^{3n-3}$ to
$(\R\times\C)^{3n-3}$.
This implies that the differential of the map
$(\boF_i^+,\boG_i^+,\boG_i^{-*},\Re(\boF_i^0),\boG_i^0)_{1\leq i\leq n-1}$ is surjective
from $(\boWp)^{3n-3}$ to $((\boWpp)^3\times\R\times\C)^{n-1}$.
Proposition \ref{monodromy-prop2} follows from the Implicit Function Theorem.\cqfd
\begin{remark}
\label{monodromy-remark1}
The kernel of the differential has real dimension $3n-3$ so we
have $3n-3$ free real parameters. These parameters correspond to deformations of
the flux vectors of the minimal $n$-noid.
\end{remark}
\section{Geometry of the immersion}
\label{geometry}
From now on, we assume that $\x(t)$ is given by Proposition \ref{monodromy-prop2}.
We write $a_{i,t}$, $b_{i,t}$ and $p_{i,t}$ for the value of the corresponding parameters.
(These parameters are in the space $\boWp$ so are functions of $\lambda$.)
For ease of notation, we write $g_t$, $\omega_t$, $\xi_t$ and $\Phi_t$ for
$g_{\x(t)}$, $\omega_{\x(t)}$, $\xi_{t,\x(t)}$ and $\Phi_{t,\x(t)}$, respectively.
Let $F_t=\Uni(\Phi_t)$.
Since the Monodromy Problem is solved, the Sym-Bobenko formula \eqref{sym-bobenko} defines
a CMC-1 immersion $f_t:\Omega\to\R^3$, where $\Omega$ is the (fixed) domain defined
by Equation \eqref{eq-Omega}.
\begin{proposition}
\label{prop-extends}
The immersion $f_t$ extends analytically to
$$\Sigma_t:=\C\cup\{\infty\}\setminus\{p_{1,t}(0),\cdots,p_{n,t}(0)\}$$
where $p_{i,t}(0)$ is the value of $p_{i,t}$ at $\lambda=0$.
\end{proposition}
We omit the proof which is exactly the same as the proof of Point 1 of Proposition 4 in \cite{nnoids}.
It relies on Theorem 3 in \cite{nnoids} which allows for $\lambda$-dependent changes of variables
in the DPW method.
\subsection{Convergence to the minimal $n$-noid}
\label{section-convergence}
\begin{proposition}
\label{convergence-prop1}
$\displaystyle\lim_{t\to 0}\smallfrac{1}{t} f_t=\psi$ where $\psi$ is (up to translation)
the conformal parametrization of the minimal $n$-noid given by Equation \eqref{minimal-eq1}.
The limit is the uniform $C^1$ convergence on compact subsets of $\Sigma_0=\C\cup\{\infty\}\setminus
\{p_{1,0},\cdots,p_{n,0}\}$.
\end{proposition}
Proof: at $t=0$, we have $g_0=g$ and $\omega_0=\omega$.
By Equation \eqref{eq-Phi0} and definition of the potential, we have
$$\Phi_0(z,\lambda)=\matrix{g(z)&1\\-1&0}\quad\mbox{ and }\quad
\frac{\partial\xi_{t;12}^{(-1)}}{\partial t}|_{t=0}=\frac{\omega}{4}.$$ 
By Theorem \ref{theorem-blowup}, $\frac{1}{t}f_t$ converges to a minimal surface with
Weierstrass data $(g,\omega)$ on compact subsets of $\Sigma_0$ minus the poles of
$g$. In a neighborhood of the poles of $g$, we use the gauge introduced in Section \ref{section-regularity}. With the notations of this section and writing $\whPhi_t=\Phi_t G_{\x(t)}$, we have
$$\whPhi_0(z,\lambda)=\matrix{1&0\\-g(z)^{-1}&1}\quad\mbox{ and }\quad
\frac{\partial\whxi_{t;12}^{(-1)}}{\partial t}|_{t=0}=\frac{g^2\omega}{4}.$$
By Theorem \ref{theorem-blowup} again, $\frac{1}{t}f_t$ converges to a minimal surface with Weierstrass data $(g,\omega)$ in a neighborhood of the poles of $g$. The two limit minimal
surfaces are of course the same, since they coincide in a neighborhood of $z_0$.\cqfd
\subsection{Delaunay ends}
\label{ends}
We denote $N_0$ the Gauss map of the minimal $n$-noid $M_0$.
For $1\leq i\leq n$, we denote $\boC_i$ the catenoid to which $M_0$ is asymptotic at $p_{i,0}$
and $\tau_i>0$ the necksize of $\boC_i$.
\begin{definition}
We say that
$N_0$ points to the inside in a neighborhood of $p_{i,0}$ if
it points to the component of $\R^3\setminus\boC_i$ containing the axis of $\boC_i$.
\end{definition}
\begin{proposition}
\label{ends-prop1}
For $1\leq i\leq n$ and $t\neq 0$:
\begin{enumerate}
\item The immersion $f_t$ has a Delaunay end at $p_{i,t}$.
If we denote $w_{i,t}$ its weight  then
$$\displaystyle\lim_{t\to 0} t^{-1}w_{i,t}=\pm 2\pi\tau_i$$
where the sign is $+$ if $N_0$ points to the inside in a neighborhood of $p_{i,0}$
and $-$ otherwise.
\item Its axis converges as $t\to 0$ to the half-line through the origin
directed by the vector $N_0(p_{i,0})$.
\item If $N_0$ points to the inside in a neighborhood of $p_{i,0}$, there exists a uniform $\varepsilon>0$ such that
for $t>0$ small enough, $f_t(D^*(p_{i,0},\varepsilon))$ is embedded.
\end{enumerate}
\end{proposition}
Proof: in  a neighborhood  of the puncture $p_{i,t}$, we may use
$w=g_t(z)-g_t(p_{i,t})$ as a local coordinate.
(This change of coordinate depends on $\lambda$. This is not a problem
by Theorem 3 in \cite{nnoids}.)
Consider the gauge
$$G(w)=\matrix{\frac{k}{\sqrt{w}} & \frac{-1}{2k\sqrt{w}}\\ 0 & \frac{\sqrt{w}}{k}}.$$
Here we can take $k=1$, but later on we will take another value of $k$ so we do the computation for general values of $k\neq 0$.
The gauged potential is
$$\whxi_t:=\xi_t\cdot G=\matrix{0 & \frac{dw}{4k^2w}+\frac{wt(\lambda-1)^2}{4k^2\lambda}\omega_t\\
\frac{k^2dw}{w} & 0}$$
Since $\omega_t$ has a double pole at $p_{i,t}$, $\whxi_t$ has a simple pole at $w=0$ with residue
$$A_{i,t}(\lambda)=\matrix{0 & \frac{1}{4k^2}+\frac{t(\lambda-1)^2}{4k^2\lambda}\alpha_{i,t}(\lambda)\\
k^2 & 0}$$
where
\begin{equation}
\label{ends-eq1}
\alpha_{i,t}=\Res_{p_{i,t}}(w\omega_t)=\Res_{p_{i,t}}(g_t(z)-g_t(p_{i,t}))\omega_t.
\end{equation}
\begin{claim}
\label{ends-claim1}
For $t$ small enough, $\alpha_{i,t}$ is a real constant (i.e. independent of $\lambda$, possibly
depending on $t$).
\end{claim}
Proof: the proof is similar to the proof of Point 2 of Proposition 4 in \cite{nnoids}.
We use the standard theory of Fuchsian systems.
Fix $t\neq 0$ and $\lambda\in\S^1\setminus\{1\}$.
Assume that $\alpha_{i,t}(\lambda)\neq 0$.
Let $\whPhi_t=\Phi_tG$.
The eigenvalues of $A_{i,t}$ are $\pm\Lambda_{i,t}$ with
$$\Lambda_{i,t}(\lambda)^2=\frac{1}{4}+\frac{t(\lambda-1)^2}{4\lambda}\alpha_{i,t}(\lambda).$$
Provided $t\neq 0$ is small enough, $\Lambda_{i,t}\not\in\Z/2$ so the system is non resonant and $\whPhi_t$ has
the following standard $z^A P$ form in the universal cover of $D(0,\varepsilon)^*$:
$$\whPhi_t(w,\lambda)=V(\lambda)\exp(A_{i,t}(\lambda)\log w)P(w,\lambda)$$
where $P(w,\lambda)$ descends to a well defined holomorphic function of $w\in D(0,\varepsilon)$ with $P(0,\lambda)=I_2$.
Consequently, its monodromy is
$$\boM_{\gamma_i}(\whPhi_t)=V(\lambda)\exp(2\pi\ii A_{i,t})V(\lambda)^{-1}$$
with eigenvalues $\exp(\pm 2\pi\ii\Lambda_{i,t}(\lambda))$.
Since the Monodromy Problem is solved, the eigenvalues are unitary complex numbers,
so $\Lambda_{i,t}(\lambda)\in\R$ which implies that $\alpha_{i,t}(\lambda)\in\R$.
This of course remains true if $\alpha_{i,t}(\lambda)=0$.
Hence $\alpha_{i,t}$ is real on $\S^1\setminus\{1\}$. 
Since all the parameters involved in the definition of $\omega_t$ are in $\boWp$,
$\alpha_{i,t}$ is holomorphic in the unit disk. Hence it is constant.\cqfd
\medskip

Returning to the proof of Proposition \ref{ends-prop1}, let $(r,s)\in\R^2$ be the solution of
\begin{equation}
\label{ends-eq2}
\left\{\begin{array}{l}
rs=\frac{1}{4}t\alpha_{i,t}\\
r+s=\frac{1}{2}\\
r<s\end{array}\right.
\end{equation}
Since $r<s$, $\sqrt{r\lambda+s}$ is well defined and does not vanish for $\lambda\in\D$.
 We take $k=\sqrt{r\lambda+s}$ in the definition of the gauge $G$. Using Equation \eqref{ends-eq2},
 we have:
$$(r\lambda^{-1}+s)(r\lambda+s)=\smallfrac{1}{4}+rs(\lambda-1)^2\lambda^{-1}=\smallfrac{1}{4}+\smallfrac{1}{4}t(\lambda-1)^2\lambda^{-1}\alpha_{i,t}.$$
So the residue of $\whxi_t$ becomes
$$A_{i,t}=
\matrix{0 & \frac{1}{r\lambda+s}\left(\frac{1}{4}+\frac{t(\lambda-1)^2}{4\lambda}\alpha_{i,t}\right)\\
r\lambda+s & 0}
=
\matrix{0 & r\lambda^{-1}+s\\ r\lambda+s & 0}$$
which is the residue of the standard Delaunay potential. By \cite{kilian-rossman-schmitt},
the immersion $f_t$ has
a Delaunay end at $p_{i,t}$ of weight
$w_{i,t}=8\pi rs=2\pi t \alpha_{i,t}$.
It remains to relate $\alpha_{i,0}$ to the logarithmic growth $\tau_i$.
For ease of notation, let us write $p_i=p_{i,0}$.
Assume that $N_0$ points to the inside in a neighborhood of $p_i$.
The flux of $M_0$ along $\gamma_i$ is equal to
$$\phi_i=2\pi\tau_i \,N_0(p_i)=2\pi\frac{\tau_i}{|g(p_i)|^2+1}\left(
2\,\Re(g(p_i)),2\,\Im(g(p_i)),|g(p_i)|^2-1\right)$$
On the other hand, we have seen in Section \ref{minimal} that the flux is equal to
$$\phi_i=-2\pi\Res_{p_i}\left(\smallfrac{1}{2}(1-g^2)\omega,\smallfrac{\ii}{2}(1+g^2)\omega,
g\omega\right)$$
Comparing these two expressions for $\phi_i$, we obtain
$$\Res_{p_i}(g\omega)=-\tau_i\frac{|g(p_i)|^2-1}{|g(p_i)|^2+1}\quad\mbox{ and }\quad
\Res_{p_i}\omega=-2\tau_i\frac{\overline{g(p_i)}}{|g(p_i)|^2+1}$$
Using Equation \eqref{ends-eq1}, this gives
$$\alpha_{i,0}=\Res_{p_i}(g\omega)-g(p_i)\Res_{p_i}\omega=\tau_i$$
If $N_0$ points to the outside in a neighborhood of $p_i$, then $\phi_i=-2\pi\tau_i\,N_0(p_i)$, so the same computation gives $\alpha_{i,0}=-\tau_i$.
This proves Point 1 of Proposition \ref{ends-prop1}.
\medskip

To prove Point 2, we use Theorem \ref{appendixB-theorem1} in Appendix
\ref{appendixB}. We need to compute $\whPhi_0$ at $w=1$. At $t=0$, we have $k=\frac{1}{\sqrt{2}}$ so
$$G(1)=\frac{1}{\sqrt{2}}\matrix{1&-1\\0&2}.$$
At $t=0$, we have $w=g(z)-g(p_i)$, so $w=1\Leftrightarrow g(z)=g(p_i)+1$.
Using Equation \eqref{eq-Phi0},
$$\whPhi_0(1)=\frac{1}{\sqrt{2}}\matrix{g(p_i)+1&1\\-1&0}\matrix{1&-1\\0&2}=\frac{1}{\sqrt{2}}\matrix{g(p_i)+1&-g(p_i)+1\\-1&1}.$$
Fix $0<\alpha<1$.
By Theorem \ref{appendixB-theorem1} (using $t\alpha_{i,t}$ as the time parameter), there exists $\varepsilon>0$, $T>0$ and $c$ such that for
$0<|t|<T$:
$$\|f_t(z)-f_{i,t}^{\boD}(z)\|\leq c|t|\,|z-p_{i,t}|^{\alpha}\quad\mbox{ in $D^*(p_{i,t},\varepsilon)$}$$
where $f_{i,t}^{\boD}:\C\setminus\{p_{i,t}\}\to \R^3$ is a Delaunay immersion.
We compute the limit axis of $f_{i,t}^{\boD}$ using Point 3 of Theorem \ref{appendixB-theorem1}:
$$\whPhi_0(1)H=\matrix{g(p_i)&1\\-1&0}=\Phi_0(p_i).$$
$$Q=F_0(p_i)$$
$$Qe_3Q^{-1}=\Nor(F_0(p_i))=N_0(p_i).$$
This proves Point 2 of Proposition \ref{ends-prop1}.
If $N_0$ points to the inside in a neighborhood of $p_{i,0}$, then for $t>0$, $t\alpha_{i,t}>0$ so Point 3 follows from Point 2 of Theorem \ref{appendixB-theorem1}.
\cqfd
\subsection{Alexandrov-embeddedness}
We recall from \cite{cosin-ros,karsten-kusner-sullivan} the definition of Alexandrov-embeddedness in the non-compact case:
\begin{definition}
\label{def-alexandrov}
A surface $M$ of finite topology is Alexandrov-embedded if M is properly immersed, if each end of $M$ is embedded, and if there exists a compact 3-manifold $\overline{W}$ with
boundary $\partial \overline{W} = \overline{S}$, $n$ points $q_1,\cdots,q_n\in\overline{S}$ and a proper immersion $F : W=\overline{W} \setminus\{q_1,\cdots, q_n\} \to\R^3$ whose restriction to $S=\overline{S}\setminus\{q_1,\cdots,q_n\}$ parametrizes $M$.
\end{definition}
Let $M$ be an Alexandrov-embedded minimal surface with $n$ catenoidal ends.
With the notations of Definition \ref{def-alexandrov},
we equip $W$ with the flat metric induced by $F$, so $F$ is a local isometry. We denote $N$ the inside normal to $S$.
\begin{lemma}
\label{lemma-alexandrov}
There exists a flat 3-manifold $W'$ containing $W$, a local isometry $F':W'\to\R^3$
extending $F$ and $r>0$ such that the tubular neighborhood $\mbox{Tub}_{r}S$
is embedded in $W'$. In other words, the map $(x,s)\mapsto \exp_x(sN(x))$ from
$S\times(-r,r)$ to $W'$ is well defined and is a diffeomorphism onto its image.
\end{lemma}
Proof:
since $M$ has catenoidal ends, there exists $r>0$ such that the inside tubular neighborhood map
$$g: S\times(0,r)\to W,\quad g(x,s)=\exp_x(s N(x))$$
is a diffeomorphism onto its image.
Since $F$ is a local isometry, we have
\begin{equation}
\label{alexandrov-eq1}
F(g(x,s))=F(x)+s\,dF(x)N(x)\quad\mbox{ for $(x,s)\in S\times(0,r)$}.
\end{equation}
We define $W'$ as the disjoint union $(S\times(-r,r))\sqcup W$ where we identify
$(x,s)\in S\times(0,r)$ with its image $g(x,s)\in W$. We define $F':W'\to\R^3$ by
$F'=F$ in $W$ and
$$F'(x,s)=F(x)+s\,dF(x)N(x)\quad\mbox{ for $(x,s)\in S\times(-r,r)$.}$$
The map $F'$ is well defined by Equation \eqref{alexandrov-eq1}.
We equip $S\times(-r,r)$ with the flat metric induced by the local diffeomorphism
$F'$, which extends the metric already defined on $S\times(0,r)$ by identification with $W$.
Since
$$dF'(x,0)(X,T)=dF(x)X+TdF(x)N(x)$$
the metric restricted to $S\times\{0\}$ is the product metric, so the normal to $S\times\{0\}$ in
$S\times(-r,r)$ is $N(x,0)=(0,1)$.
Since $F'$ is a local isometry, we have for $(x,s)\in S\times (-r,r)$
$$F'\left(\exp_{(x,0)}sN(x,0)\right)=F'(x,0)+sdF'(x,0)(0,1)=F(x)+sdF(x)N(x)=F'(x,s)$$
Hence $\exp_{(x,0)}sN(x,0)=(x,s)$ so
$\mbox{Tub}_{r}(S\times\{0\})$ is embedded in $S\times(-r,r)$.\cqfd
\medskip

\begin{figure}
\begin{center}
\includegraphics[height=3cm]{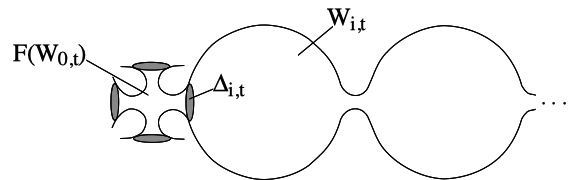}
\end{center}
\caption{Decomposition of a 4-noid into pieces. Only one Delaunay end is represented, and $F(W_{0,t})$ is represented as an embedded domain for clarity, but in general it will be immersed.}
\label{fig2}
\end{figure}

We now return to the proof of Theorem \ref{intro-theorem1}.
We orient the minimal $n$-noid $M_0$ so that its Gauss map points to the inside in a neighborhood of $p_1$.
For $0<|t|<\epsilon$, we denote $M_t$ the image of the immersion $f_t$ that we have constructed.
\begin{proposition}
\label{prop-alexandrov} If $M_0$ is Alexandrov embedded, then
for $t>0$ small enough, $M_t$ is Alexandrov embedded.
\end{proposition}
Proof: our strategy is to cut $M_t$ by suitable planes into pieces which are either close to $M_0$ or Delaunay surfaces (see Figure \ref{fig2}). Then we prove that each piece, together with flat disks in the cutting planes, is the boundary of a domain, using the Jordan Brouwer Theorem.
\medskip

Since $M_0$ is Alexandrov embedded,
$N_0$ points to the inside in a neighborhood of each end, so $M_t$ has embedded ends
by Proposition \ref{ends-prop1}.
Let $\varepsilon>0$ be the number given by our application of Theorem \ref{appendixB-theorem1}
in Section \ref{ends} and $f_{i,t}^{\boD}:\C\setminus\{p_{i,t}\}\to \R^3$ be the Delaunay immersion
which approximates $f_t$ in $D^*(p_{i,t},\varepsilon)$.
Recall that $f_t(D^*(p_{i,t},\varepsilon))$ is embedded.
Let $\wtf_t=\frac{1}{t}f_t$. By Proposition \ref{convergence-prop1},
$\wtf_t$ converges to $\psi$ on compact subsets of $\Sigma_0$, where
$\psi:\Sigma_0\to\R^3$ is a parametrization of $M_0$.
Since $M_0$ has catenoidal ends, we may assume (taking $\varepsilon$ smaller if necessary) that $\psi(D^*(p_{i,0},\varepsilon))$ is embedded and $N_0\neq N_0(p_{i,0})$ in
$D^*(p_{i,0},\varepsilon)$.
\medskip

Let $h_i:\R^3\to\R$ be the height function in the direction $N_0(p_{i,0})$, defined by
$$h_i(x)=\langle x,N_0(p_{i,0})\rangle.$$
We shall cut $M_t$ by the plane $h_i=\delta$ where
$\delta>0$ is a fixed, large enough number such that for $1\leq i\leq n$,
$$\delta>\max_{C(p_{i,0},\varepsilon)}h_i\circ\psi.$$
Since $\displaystyle\lim_{z\to p_{i,0}} h_i\circ \psi(z)=+\infty$, we may fix a positive, small enough $\varepsilon'<\varepsilon$ such that
$$\min_{C(p_{i,0},\varepsilon')}h_i\circ\psi>\delta.$$
Let $\boA_{i,t}$ be the annulus defined by $\varepsilon'\leq |z-p_{i,t}|\leq \varepsilon$.
Since $N_0\neq N_0(p_{i,0})$ in $\boA_{i,0}$,
$$\min_{\boA_{i,0}} \|N_0(z)-N_0(p_{i,0})\|>0.$$
For $t>0$ small enough:
\begin{equation}
\label{alexandrov-eq3}
\max_{C(p_{i,t},\varepsilon)}h_i\circ\wtf_t<\delta
\end{equation}
\begin{equation}
\label{alexandrov-eq4}
\min_{C(p_{i,t},\varepsilon')}h_i\circ\wtf_t>\delta
\end{equation}
\begin{equation}
\label{alexandrov-eq5}
\min_{\boA_{i,t}} \|N_t(z)-N_0(p_{i,0})\|>0.
\end{equation}
Hence the function $h_i\circ\wtf_t$ has no critical point in the annulus $\boA_{i,t}$.
So $h_i\circ \wtf_t=\delta$ defines a regular closed curve $\gamma_{i,t}$ in $\boA_{i,t}$. At $t=0$,
$h_i\circ\psi=\delta$ is a single curve around $p_{i,0}$, so $\gamma_{i,t}$ has only one component and is not contractible in $\boA_{i,t}$. Let $D_{i,t}\subset\C$ be the topological disk bounded by $\gamma_{i,t}$
and $D_{i,t}^*=D_{i,t}\setminus\{p_{i,t}\}$.
Let $\Delta_{i,t}$ be the closed topological disk bounded by $\wtf_t(\gamma_{i,t})$
in the plane defined by $h_i(x)=\delta$.
\begin{claim}
\label{claim1}
For $t>0$ small enough, $\wtf_t(D_{i,t}^*)\cap \Delta_{i,t}=\emptyset$.
\end{claim}
Proof: of course, $h_i\circ \wtf_t>\delta$ in $D_{i,t}^*\cap\boA_{i,t}$. What we need to prove is that
$\wtf_t(D^*(p_{i,t},\varepsilon'))$ does not intersect $\Delta_{i,t}$. We do this by comparison with
the Delaunay surface.
Let $\Pi_i=N_0(p_{i,0})^{\perp}$ and $\pi_i=\R^3\to \Pi_i$ be the orthogonal projection.
Since $\psi$ has a catenoidal end at $p_{i,0}$, $\psi(\boA_{i,t})$ is a graph over an annulus in the plane $\Pi_i$, with inside boundary circle $\pi_i\circ\psi(C(p_{i,t},\varepsilon))$ and outside boundary circle
$\pi_i\circ\psi(C(p_{i,t},\varepsilon'))$. Moreover, $N_0$ is close to $N_0(p_{i,t})$.
Since $\wtf_t$ is $C^1$ close to $\psi$ in $\boA_{i,t}$, for $t>0$ small enough,
$\wtf_t(\boA_{i,t})$ is a graph over
an annulus in the plane $\Pi_i$, with inside boundary circle $\pi_i\circ\wtf_t(C(p_{i,t},\varepsilon))$ and
outside boundary circle $\pi_i\circ\wtf_t(C(p_{i,t},\varepsilon'))$.
\medskip

Now we go back to the original scale.
Since $f_t$ is $C^1$ close to $f_{i,t}^{\boD}$ in $D^*(p_{i,t},\varepsilon)$, we conclude that
$f_{i,t}^{\boD}(\boA_{i,t})$ is a graph over an annulus in the plane $\Pi_i$, with inside boundary circle
$\pi_i\circ f_{i,t}^{\boD}(C(p_{i,t},\varepsilon))$ and outside boundary circle $\pi_i\circ f_{i,t}^{\boD}(C(p_{i,t},\varepsilon'))$.
Then from the geometry of Delaunay surfaces, there exists a curve
$\gamma_{i,t,0}$ in $D^*(p_{i,t},\varepsilon')$ such that $f_{i,t}^{\boD}(\gamma_{i,t,0})$ is
a closed curve in the plane $h_i=\frac{1}{2}$. Let $D_{i,t,0}$ be the disk bounded by $\gamma_{i,t,0}$ and $\boA_{i,t,0}$ be the closed annulus bounded by $\gamma_{i,t}$ and $\gamma_{i,t,0}$.
Then $h_i\circ f_{i,t}^{\boD}>\frac{1}{2}$ in $D_{i,t,0}$ and $f_{i,t}^{\boD}(\boA_{i,t,0})$ is a graph over an annulus in the plane $\Pi_i$.
Since $f_t$ is $C^1$ close to $f_{i,t}^{\boD}$ in $D^*(p_{i,t},\varepsilon)$, we conclude that
$h_i\circ f_t>\frac{1}{4}$ in $D_{i,t,0}^*$ and $f_t(\boA_{i,t,0})$ is a graph over an annulus in the plane $\Pi_i$.

\medskip
Back to the scale $\frac{1}{t}$, $\wtf_t(D_{i,t}\cap\boA_{i,t,0})$ is a graph over an annulus in the plane
$\Pi_i$ whose inside boundary circle is $\pi_i\circ\wtf_t(\gamma_{i,t})=\pi_i(\partial\Delta_{i,t})$, so
$\wtf_t(D_{i,t}\cap\boA_{i,t,0})\cap\Delta_{i,t}=\emptyset$. Moreover,
$h_i\circ\wtf_t>\frac{1}{4t}\gg \delta$ in $D_{i,t,0}^*$ so $\wtf_t(D^*_{i,t,0})\cap\Delta_{i,t}=\emptyset$.
\cqfd
\medskip

\begin{claim}
For $t>0$ small enough,
$\wtf_t(D_{i,t}^*)\cup \Delta_{i,t}$ is the boundary of a cylindrically bounded
domain $W_{i,t}\subset\R^3$.
\end{claim}
Proof: since $f_t$ is close to $f_{i,t}^{\boD}$ in $D_{i,t}^*$, we can
find an increasing diverging sequence $(R_{t,k})_{k\in\N}$ such that $f_t(D_{i,t}^*)$ intersects the plane $h_i=R_{t,k}$
transversally along a closed curve $f_t(\gamma_{i,t,k})$.
(Explicitely, we can take $R_{t,k}=\frac{1}{2}+k\,h_i(T_t)$ where $T_t\in\R^3$ is the period of
the Delaunay surface $f_{i,t}^{\boD}$.)
Let $\boA_{i,t,k}$ be the annulus bounded by
$\gamma_{i,t}$ and $\gamma_{i,t,k}$. Let $\Delta_{i,t,k}$ be the closed disk
bounded by $\wtf_t(\gamma_{i,t,k})$ in the plane $h_i=t^{-1}R_{t,k}$.
Then $\wtf_t(\boA_{i,t,k})\cup\Delta_{i,t}\cup\Delta_{i,t,k}$ is topologically a sphere: the image of $\S^2$ by an injective continuous map. By the Jordan Brouwer Theorem, it is the boundary of a bounded domain $W_{i,t,k}$. Clearly, $W_{i,t,k}\subset W_{i,t,k+1}$.
We take $W_{i,t}=\bigcup_{k\in\N} W_{i,t,k}$.
\cqfd
\medskip

Let $\Omega_t=\C\cup\{\infty\}\setminus (D_{1,t}\cup\cdots\cup D_{n,t})$.
Let $W'$ be the flat 3-manifold given by Lemma \ref{lemma-alexandrov}
and denote $F:W'\to\R^3$ its developing map (instead of $F'$).
(Here $W'$ is an open manifold, meaning not a manifold-with-boundary.)
\begin{claim}
For $t>0$ small enough, there exists a compact domain $W_{0,t}$ in $W'$ such that
$$F(\partial W_{0,t})=\wtf_t(\Omega_t)\cup\Delta_{1,t}\cup\cdots\cup\Delta_{n,t}.$$
\end{claim}
Proof: by definition, $\psi$ lifts to a diffeomorphism $\whpsi:\Sigma_0\to S\subset W'$
such that $F\circ \whpsi=\psi$.
Since $M_0$ has catenoidal ends, there exists domains $V_1,\cdots,V_n$ in $W'$ such that
 for $1\leq i\leq n$:
\begin{itemize}
\item $F:V_i\to F(V_i)\subset\R^3$ is a diffeomorphism,
\item $V_i$ is foliated by flat disks on which $h_i\circ F$ is constant (in particular, $h_i\circ F$ is constant on $\partial V_i$),
\item $\whpsi(D^*(p_{i,0},\varepsilon))\subset V_i$ (which might require taking a
smaller $\varepsilon>0$),
\item $h_i<\delta$ on $V_i\cap \whpsi(\Sigma_0\setminus\bigcup_{i=1}^n D(p_{i,0},\varepsilon))$
(which might require taking a larger $\delta$).
\end{itemize}
Let $r>0$ be the radius of the embedded tubular
neighborhood of $S$ in $W'$ constructed in Lemma \ref{lemma-alexandrov}.
For $t>0$ small enough, $||\wtf_t-\psi||<r$ in $\overline{\Omega}_t$, so $\wtf_t$ lifts to
$\whf_t:\overline{\Omega}_t\to W'$ such that $F\circ\whf_t=\wtf_t$.
(Explicitely, $\whf_t(z)=\exp_{\psi(z)}(\wtf_t(z)-\psi(z)).$)
From the properties of $V_i$ and the convergence of $\whf_t$ to $\whpsi$ on compact
subsets of $\Sigma_0$, we have for $t>0$ small enough
\begin{equation}
\label{alexandrov-eq6}\whf_t(\overline{\Omega}_t\cap D(p_{i,0},\varepsilon))\subset V_i
\end{equation}
\begin{equation}
\label{alexandrov-eq7}
h_i<\delta \quad\mbox{ on $V_i\cap \whf_t(\Sigma_0\setminus\bigcup_{i=1}^n D(p_{i,0},\varepsilon))$}.
\end{equation}
By Equation \eqref{alexandrov-eq6}, $\whf_t(\gamma_{i,t})\subset V_i$ so
$\Delta_{i,t}$ lifts to a closed disk $\whDelta_{i,t}\subset V_i$ such that $\partial\whDelta_{i,t}=\whf_t(\gamma_{i,t})$ and $F(\whDelta_{i,t})=\Delta_{i,t}$.
Since $F$ is a diffeomorphism on $V_i$,
$\whf_t(\Omega_t\cap D(p_{i,t},\varepsilon))$ is disjoint from $\whDelta_{i,t}$.
By \eqref{alexandrov-eq7}, $\whf_t(\Omega_t\setminus\bigcup_{i=1}^n D(p_{i,t},\varepsilon))$ is disjoint from
$\whDelta_{i,t}$. Hence $\whf_t(\Omega_t)\cap\whDelta_{i,t}=\emptyset$.
Then $\whf_t(\Omega_t)\cup\whDelta_{1,t}\cup\cdots\cup\whDelta_{n,t}$ is a topological sphere
in $W'$.
Since $M_0$ has genus zero, $W'$ is homeomorphic to $\R^3$. By the Jordan Brouwer Theorem, $\whf_t(\Omega_t)\cup\whDelta_{1,t}\cup\cdots\cup\whDelta_{n,t}$ is the boundary of a compact domain $W_{0,t}\subset W'$.
\cqfd
\medskip

Returning to the proof of Proposition \ref{prop-alexandrov}, let $W_t$ be the abstract 3-manifold with boundary obtained as the disjoint union $\overline{W}_{0,t}\sqcup \overline{W}_{1,t}\sqcup\cdots\sqcup \overline{W}_{n,t}$,
identifying $\overline{W}_{0,t}$ and $\overline{W}_{i,t}$ along their boundaries $\whDelta_{i,t}$ and
$\Delta_{i,t}$ via the map $F$ for $1\leq i\leq n$.
Let $F_t:W_t\to\R^3$ be the map defined by $F_t=F$ in $\overline{W}_{0,t}$ and $F_t=\mbox{id}$ in
$\overline{W}_{i,t}$ for $1\leq i\leq n$. Then $F_t$ is a proper local diffeomorphism
whose boundary restriction parametrizes $M_t$.
Moreover, since each $\overline{W}_{i,t}$ is homeomorphic to a closed ball minus a boundary point,
we may compactify $W_t$ by adding $n$ points. This proves that $M_t$ is Alexandrov-embedded.
\cqfd
\appendix
\section{Appendix: complements on the Banach algebra $\boW$}
\label{appendix-iwasawa}
In this section, we prove several basic facts about the Banach algebra $\boW$
introduced in Section \ref{section-functional-spaces} that are used in this paper and related papers
\cite{nnoids,traizet3}.
\begin{proposition}
\label{appendix-prop-algebra}
If $\rho\geq 1$, $\boW$ is a Banach algebra.
\end{proposition}
Proof:
let $f,g\in\boW$ and $h=fg$. Then $(|f_i|\rho^{|i|})_{i\in\Z}$ and
$(|g_i|\rho^{|i|})_{i\in\Z}$ are summable families, so
$(|f_ig_j|\rho^{|i|+|j|})_{(i,j)\in\Z^2}$ is a summable family.
Using the triangular inequality and $\rho\geq 1$, we obtain that
$|f_ig_j|\rho^{|i+j|}$ is a summable family, so the following computation is valid:
$$\sum_{i\in\Z}|h_i|\rho^{|i|}
\leq\sum_{(i,k)\in\Z^2}|f_kg_{i-k}|\rho^{|i|}
=\sum_{(i,j)\in\Z^2}|f_ig_j|\rho^{|i+j|}
\leq\sum_{i\in\Z}|f_i|\rho^{|i|}\sum_{j\in\Z}|g_j|\rho^{|j|}.$$
Hence $h\in\boW$ and $\|h\|\leq \|f\|\times\|g\|$.\cqfd
\medskip

Recall that if $L$ is a loop group, $L_{\rho}\subset L$ denotes the subgroup of loops whose
entries are in the Banach algebra $\boW$.
\begin{proposition}
\label{appendix-iwasawa-proposition1}
Let $\Phi\in\Lambda SL(2,\C)_{\rho}$ and $(F,B)$ be its Iwasawa decomposition.
Then $F\in\Lambda SU(2)_{\rho}$ and $B\in\Lambda_+^{\R}SL(2,\C)_{\rho}$.
\end{proposition}
Proof: $\Phi$ extends holomorphically to the annulus $\frac{1}{\rho}<|\lambda|<\rho$
and $B$ extends holomorphically to the disk $|\lambda|<1$, so
$F=\Phi B^{-1}$ extends holomorphically to the annulus
$\frac{1}{\rho}<|\lambda|<1$.
By an application of the Schwarz reflection principle, the fact that $F(\lambda)\in SU(2)$
on the unit circle implies that $F$ extends holomorphically to the annulus
$\frac{1}{\rho}<|\lambda|<\rho$ and satisfies (see details in Appendix A of \cite{raujouan}):
$$F_{11}(\smallfrac{1}{\overline{\lambda}})=\overline{F_{22}(\lambda)}\quad\mbox{ and }\quad
F_{12}(\smallfrac{1}{\overline{\lambda}})=-\overline{F_{21}(\lambda)}.$$
We expand $F=(F_{k\ell})_{1\leq k,\ell\leq 2}$ in the annulus $\frac{1}{\rho}<|\lambda|<\rho$ as
$$F(\lambda)=\sum_{i\in\Z} F_i\lambda^i\quad\mbox{ with }\quad
F_i=(F_{k\ell,i})_{1\leq k,\ell\leq 2}.$$
Then for all $i\in\Z$
\begin{equation}
\label{appendixA-eqSU2}
F_{11,-i}=\overline{F_{22,i}}\quad\mbox{ and }\quad
F_{12,-i}=-\overline{F_{21,i}}.
\end{equation}
We expand $\Phi$ and $B^{-1}$ in the annulus $\frac{1}{\rho}<|\lambda|<1$ as
$$\Phi(\lambda)=\sum_{i\in\Z} \Phi_i \lambda^i,\qquad
B^{-1}(\lambda)=\sum_{i\in\Z} C_i\lambda^i$$
with $C_i=0$ for $i<0$.
Then
$$F_i=\sum_{j\in\Z} \Phi_j C_{i-j}.$$
Since $\Phi\in\Lambda SL(2,\C)_{\rho}$,
$(\|\Phi_i\|\rho^{-i})_{i\in\Z}$ is a summable family
(here $\|A\|=\max |A_{k\ell}|$ for $A\in\boM_2(\C)$).
Since $B^{-1}$ is holomorphic in the unit disk,
$(\|C_i\|\rho^{-i})_{i\in\Z}$ is a summable family.
Hence
$$\sum_{i\in\Z} \|F_i\|\rho^{-i}<\infty.$$
By Equation \eqref{appendixA-eqSU2}, this implies
$$\sum_{i\in\Z}\|F_i\|\rho^i<\infty.$$
Hence $F\in\boM_2(\boW)$.
Since $\boW$ is a Banach algebra, $B=F^{-1}\Phi\in\boM_2(\boW)$ as well.\cqfd
\begin{remark}
If $\Phi\in\Lambda SL(2,\C)$ is a Laurent polynomial of degree $d$,
a similar argument proves that $F$ is a Laurent polynomial of degree at most $d$ and
$B$ is a polynomial of degree at most $2d$.
\end{remark}
The loop groups $\Lambda SL(2,\C)_{\rho}$, $\Lambda SU(2)_{\rho}$ and
$\Lambda_+^{\R} SL(2,\C)_{\rho}$ are Banach manifolds. This can be proved using the submersion criterion
(see \cite{bourbaki}, 5.9.1 for the definition of a submersion in the Banach case).
The tangent space at $I_2$ of these Banach Lie groups are respectively the following Banach Lie algebras:
$$\Lambda\sl(2,\C)_{\rho}=\{M\in\boM_2(\boW):M_{11}+M_{22}=0\}$$
$$\Lambda\su(2)_{\rho}=\{M\in\boM_2(\boW):M_{11}+M_{22}=0,\;M_{11}+M_{11}^*=0,\;M_{12}+M_{21}^*=0\}$$
$$\Lambda_+^{\R}\sl(2,\C)_{\rho}=\{M\in\boM_2(\boWp):M_{11}+M_{22}=0,\;
M_{11}(0)\in\R,\;M_{21}(0)=0\}.$$
\begin{proposition}
\label{appendix-iwasawa-proposition2}
$$\Lambda \sl(2,\C)_{\rho}=\Lambda\su(2)_{\rho}\oplus\Lambda_+^{\R}\sl(2,\C)_{\rho}.$$
\end{proposition}
Proof:
\begin{itemize}
\item Let $M\in\Lambda\su(2)_{\rho}\cap\Lambda_+^{\R}\sl(2,\C)_{\rho}$. Then all entries of $M$ are in $\boWp\cap\boWn=\boWo$ so $M$ is constant. It is then straightforward that $M=0$.
\item Let $M\in\Lambda\sl(2,\C)_{\rho}$ and write
$$M(\lambda)=\sum_{n\in\Z}M_n\lambda^n\quad\mbox{ with }\quad  M_0=\matrix{a_0&b_0\\c_0&-a_0}.$$
Define
$$\duni(M)(\lambda)=\matrix{\ii\,\Im(a_0) & -\overline{c_0}\\c_0 & -\ii\,\Im(a_0)}
+\sum_{n<0}\left(M_n\lambda^n-\mbox{}^t \overline{M}_n\lambda^{-n}\right)\in \Lambda\su(2)_{\rho}$$
$$\dpos(M)(\lambda)=\matrix{\Re(a_0)&b_0+\overline{c_0}\\ 0 & -\Re(a_0)}
+\sum_{n>0}M_n\lambda^n+\sum_{n<0}\mbox{}^t \overline{M}_n\lambda^{-n}\in\Lambda_+^{\R}\sl(2,\C)_{\rho}.$$
Then
$M=\duni(M)+\dpos(M)$.\cqfd
\end{itemize}
\begin{theorem}
\label{appendix-iwasawa-theorem}
Iwasawa decomposition is a smooth diffeomorphism (in the sense of Banach manifolds) from
$\Lambda SL(2,\C)_{\rho}$ to
$\Lambda SU(2)_{\rho}\times\Lambda_+^{\R}SL(2,\C)_{\rho}$.
Moreover, its differential at $I_2$ is given by
$d\Uni(I_2)=\duni$ and $d\Pos(I_2)=\dpos$.
\end{theorem}
Proof: consider the following smooth map:
$$p:\Lambda SU(2)_{\rho}\times\Lambda_+^{\R}SL(2,\C)_{\rho}\to\Lambda SL(2,\C)_{\rho},\quad p(F,B)=FB.$$
Then
$$dp(I_2,I_2):\Lambda\su(2)_{\rho}\times\Lambda_+^{\R}\sl(2,\C)_{\rho}\to\Lambda\sl(2,\C)_{\rho},\quad dp(I_2,I_2)(F,B)=F+B.$$
By Proposition \ref{appendix-iwasawa-proposition2}, $dp(I_2,I_2)$ is an isomorphism.
By the Inverse Mapping Theorem, $p$ is a local diffeomorphism in a neighborhood of
$(I_2,I_2)$.
Let $(F_0,B_0)\in\Lambda SU(2)_{\rho}\times\Lambda_+^{\R}SL(2,\C)_{\rho}$.
Consider the maps
$$f_1:\Lambda SU(2)_{\rho}\times\Lambda_+^{\R}SL(2,\C)_{\rho}\to\Lambda SU(2)_{\rho}\times\Lambda_+^{\R}SL(2,\C)_{\rho},\quad
f_1(F,B)=(F_0F,BB_0)$$
$$f_2:\Lambda SL(2,\C)_{\rho}\to\Lambda SL(2,\C)_{\rho},\quad f_2(\Phi)=F_0\Phi B_0.$$
Then $f_1$ is a diffeomorphism with inverse $(F,B)\mapsto (F_0^{-1}F,BB_0^{-1})$
and $f_2$ is a diffeomorphism with inverse $\Phi\mapsto F_0^{-1}\Phi B_0^{-1}$.
We have
$$f_2\circ p\circ f_1^{-1}=p$$
so $p$ is a local diffeomorphism in a neighborhood of any $(F_0,B_0)$.
Since $p$ is injective (by the standard Iwasawa decomposition theorem)
and onto (by Proposition \ref{appendix-iwasawa-proposition1}),
$p$ is a smooth diffeomorphism. Iwasawa decomposition is of course the inverse of $p$.\cqfd
\medskip

As an application of Theorem \ref{appendix-iwasawa-theorem}, let us prove Equation
\eqref{eq-df} by differentiation of the Sym-Bobenko formula \eqref{sym-bobenko}.
In the following computation, we omit the $\lambda$ variable.
Fix $z_0\in\Sigma$ and define
$\wtPhi(z)=F(z_0)^{-1}\Phi(z)B(z_0)^{-1}$.
Then
$$d\wtPhi(z_0)=F(z_0)^{-1}\Phi(z_0)\xi(z_0)B(z_0)^{-1}=B(z_0)\xi(z_0)B(z_0)^{-1}.$$
Writing $B(z_0)^0=\minimatrix{\rho&\mu\\0&\rho^{-1}}$ and $\xi_{12}(z_0)=\lambda^{-1}\beta$, we have
$$d\wtPhi(z_0)^-=\matrix{\rho&\mu\\0&\rho^{-1}}\matrix{0&\lambda^{-1}\beta^0\\0&0}\matrix{\rho^{-1}&-\mu\\0&\rho}=\matrix{0&\rho^2\lambda^{-1}\beta^0\\0&0}.$$
Since $\wtPhi(z_0)=I_2$, Theorem \ref{appendix-iwasawa-theorem} yields
$$F(z_0)^{-1}dF(z_0)=d\wtF(z_0)=\duni(d\wtPhi(z_0))=\matrix{0&\rho^2\lambda^{-1}\beta^0\\-\rho^2\lambda\overline{\beta^0}&0}+\mbox{constant}$$
$$\frac{\partial}{\partial\lambda}\left(F(z_0)^{-1}dF(z_0)\right)|_{\lambda=1}=\matrix{0&-\rho^2\beta^0\\-\rho^2\overline{\beta^0}&0}.$$
By differentiation of the Sym-Bobenko formula, we obtain (omitting $z_0$)
\begin{eqnarray*}
df&=&-2\ii\left(\frac{\partial dF}{\partial\lambda}F^{-1}-\frac{\partial F}{\partial\lambda}F^{-1}dF F^{-1}\right)|_{\lambda=1}
\\
&=&-2\ii F\frac{\partial}{\partial\lambda}\left(F^{-1}dF\right)F^{-1}|_{\lambda=1}
\\
&=&2\ii \rho^2 F\matrix{0&\beta^0\\ \overline{\beta^0}&0}F^{-1}|_{\lambda=1}.
\end{eqnarray*}
\section{Appendix: On Delaunay ends in the DPW method}
\label{appendixB}
We consider the standard Delaunay residue for $t\leq \frac{1}{16}$:
$$A_t(\lambda)=\matrix{0&\lambda^{-1}r+s\\\lambda r+s&0}\qquad
\mbox{ where }
\left\{\begin{array}{l}
r+s=\frac{1}{2}\\
rs=t\\
r<s\end{array}\right.$$
In particular, in the limit case $t=0$, we have
$$A_0=\matrix{0&\frac{1}{2}\\\frac{1}{2}&0}.$$
Let $\A_{\rho}$ be the annulus $\frac{1}{\rho}<|\lambda|<\rho$, where $\rho>1$.
\begin{definition}[\cite{raujouan}]
A perturbed Delaunay potential is a family of DPW potentials $\xi_t$ of the form
$$\xi_t(z,\lambda)=A_t(\lambda)\frac{dz}{z}+R_t(z,\lambda)dz$$
where $R_t$ is of class $C^2$ with respect to $(t,z,\lambda)\in (-T,T)\times D(0,\varepsilon)\times\A_{\rho}$ for some positive $\varepsilon$ and $T$, and satisfies $R_0=0$. In particular, $\xi_0=A_0\frac{dz}{z}$.
\end{definition}
Let $(e_1,e_2,e_3)$ represent the canonical basis of $\R^3$ in the $\su(2)$ model.
\begin{theorem}
\label{appendixB-theorem1}
Let $\xi_t$ be a perturbed Delaunay potential.
Let $\Phi_t(z,\lambda)$ be a family of solutions of $d\Phi_t=\Phi_t\xi_t$ in the universal cover of the punctured
disk $D^*(0,\varepsilon)$.
Assume that $\Phi_t(z,\lambda)$ depends continuously on $(t,z,\lambda)$ and that
the Monodromy Problem for $\Phi_t$ is solved.
Let $f_t=\Sym(\Uni(\Phi_t))$ be the immersion given by the DPW method.
Finally, assume that $\Phi_0(1,\cdot)$ is constant (i.e. independent of $\lambda$).
\medskip

Given $0<\alpha<1$, there exists uniform positive numbers $\varepsilon'\leq \varepsilon$, $T'\leq T$, $c$ and a family of Delaunay immersions
$f_t^{\boD}:\C^*\to\R^3$ such that:
\begin{enumerate}
\item For $0<|t|<T'$ and $0<|z|<\varepsilon'$:
$$\|f_t(z)-f_t^{\boD}(z)\|\leq c|t|\,|z|^\alpha.$$
\item For $0<t<T'$, $f_t:D^*(0,\varepsilon')\to\R^3$ is an embedding.
\item The end of $f_t^{\boD}$ at $z=0$ has weight $8\pi t$ and its axis converges when $t\to 0$
to the half-line spanned by the vector $Qe_3Q^{-1}$ where
$$Q=\Uni\left(\Phi_0(1)H\right)
\quad\mbox{ and }\quad
H=\frac{1}{\sqrt{2}}\matrix{1&1\\-1&1}.$$
\end{enumerate}
\end{theorem}
Thomas Raujouan has proved this result in \cite{raujouan}, Theorem 3, in the case
$\Phi_0(1,\lambda)=I_2$. He proves that the limit axis is spanned by
$e_1$.
(In fact, he finds that the limit axis is $-e_1$, but this is because he has the opposite sign in the Sym-Bobenko formula. See Remark \ref{remark-signs}.)
Then in Section 2 of \cite{raujouan}, he explains, in the case $r>s$, how to extend his result to
the case where $\Phi_0(1,\lambda)$ is constant. We adapt his method to the case $r<s$.
\begin{lemma}
\label{appendixB-lemma1}
There exists a gauge $G(z)$ and a change of variable $h(z)$ with $h(0)=0$ 
such that $\wtxi_t=(h^*\xi_t)\cdot G$ is a perturbed Delaunay potential (with residue $A_t$) 
and $\wtPhi_t=(h^*\Phi_t)\times G$
satisfies at $t=0$
\begin{equation}
\label{appendixB-eq2}
\wtPhi_0(1,\lambda)=QH^{-1}\in\Lambda SU(2).\end{equation}
\end{lemma}
Proof: we follow the method explained in Section 2 of \cite{raujouan}.
We take the change of variable in the form
$$h(z)=\frac{z}{pz+q}$$
where $p,q$ are complex numbers (independent of $t$) to be determined, with $q\neq 0$.
We consider the following gauge:
$$G(z)=\frac{1}{\sqrt{q(pz+q)}}\matrix{pz+q&pz\\0&q}.$$
It is chosen so that
\begin{equation}
\label{appendixB-eq3}
G(0)=I_2\quad\mbox{ and }\quad dG=GA_0\frac{dz}{z}-A_0 G\frac{dh}{h}.\end{equation}
(In fact, the gauge $G$ is found as the only solution of Problem \eqref{appendixB-eq3} which is upper triangular.)
We have
$$\wtxi_t=(h^*\xi_t)\cdot G=G^{-1}\left(A_t(\lambda)\frac{dh}{h}+R_t(h,\lambda)dh\right)G+G^{-1}dG.$$
Since $G(0)=I_2$, $\wtxi_t$ has a simple pole at $z=0$ with residue
$A_t$.
Using Equation \eqref{appendixB-eq3}, we obtain at $t=0$:
$$\wtxi_0=G^{-1}A_0\frac{dh}{h}G+G^{-1}dG=A_0\frac{dz}{z}.$$
Hence $\wtxi_t$ is a perturbed Delaunay potential.
It remains to compute $\wtPhi_0(1)$.
The matrix $H$ diagonalises $A_0$:
$$A_0=H\matrix{\frac{-1}{2}&0\\0&\frac{1}{2}}H^{-1}$$
Hence
$$\Phi_0(z)=\Phi_0(1)z^{A_0}=\Phi_0(1)H\matrix{\frac{1}{\sqrt{z}}&0\\0&\sqrt{z}}H^{-1}$$
\begin{eqnarray*}
\wtPhi_0(1)&=&\Phi_0(h(1))G(1)\\
&=&\Phi_0(1)H\matrix{\sqrt{p+q}&0\\0&\frac{1}{\sqrt{p+q}}}
H^{-1}\frac{1}{\sqrt{q(p+q)}}\matrix{p+q&p\\0&q}\\
&=&\Phi_0(1)H\matrix{\sqrt{q}&\frac{p}{\sqrt{q}}\\0&\frac{1}{\sqrt{q}}}H^{-1}
\end{eqnarray*}
We decompose
$\Phi_0(1)H=QR$ with $Q\in SU(2)$ and 
$R=\minimatrix{\rho&\mu\\0&\frac{1}{\rho}}$.
Then
$$\wtPhi_0(1)=Q\matrix{\rho&\mu\\0&\frac{1}{\rho}}\matrix{\sqrt{q}&\frac{p}{\sqrt{q}}\\0&\frac{1}{\sqrt{q}}}H^{-1}$$
We take $q=\frac{1}{\rho^2}$ and $p=-\frac{\mu}{\rho}$ to cancel the two matrices in the middle and obtain
Equation \eqref{appendixB-eq2}.
\cqfd
\medskip

We can now prove Theorem \ref{appendixB-theorem1}.
Let
$$\whPhi_t(z,\lambda)=HQ^{-1}\wtPhi_t(z,\lambda)=HQ^{-1}\Phi_t(h(z),\lambda)G(z,\lambda)$$
Since $\whPhi_0(1,\lambda)=I_2$,
we can apply Theorem 3 in \cite{raujouan} which says that the resulting immersion $\widehat{f}_t$ satisfies Points 1 and 2 of Theorem \ref{appendixB-theorem1} and its limit axis is
spanned by $e_1$. We have
$$f_t\circ h=QH^{-1}\widehat{f_t} HQ^{-1}$$
so $f_t\circ h$ and $\widehat{f}_t$ differ by a rotation and the limit axis of $f_t$ is spanned by
the vector $QH^{-1}e_1HQ^{-1}$.
Now
$$H^{-1}e_1H=-\frac{\ii}{2}\matrix{1&-1\\1&1}\matrix{0&1\\1&0}\matrix{1&1\\-1&1}
=-\ii \matrix{-1&0\\0&1}=e_3.$$
\cqfd

\noindent
Martin Traizet\\
Institut Denis Poisson\\
Universit\'e de Tours, 37200 Tours, France\\
\verb$martin.traizet@univ-tours.fr$
\end{document}